\documentclass[11pt]{article}
\usepackage{amsmath} % Provides align environment
\usepackage{bm}
\usepackage{authblk}
\usepackage{array}
\usepackage{booktabs}
\usepackage{graphicx} 
\usepackage{float} % for H placement specifier
\usepackage{tabularx}
\newcolumntype{L}{>{\raggedright\let\newline\\\arraybackslash\hspace{0pt}}X}
\usepackage[margin=1in]{geometry} % Adjust margins here
\usepackage{changepage}           % for the Key Points spacings
\usepackage{setspace}             % for onehalfspacing
\usepackage{caption}
\usepackage[T1]{fontenc}   % better font encoding for accented characters
\usepackage[utf8]{inputenc} % allows UTF-8 characters

\usepackage[round]{natbib}  % or [numbers] if you prefer numeric citations

\usepackage[dvipsnames]{xcolor}
\usepackage[colorlinks]{hyperref}
\usepackage{siunitx}
\usepackage[textsize=footnotesize]{todonotes}
\usepackage[acronym]{glossaries}
\glsdisablehyper  % don't use a hyperlink for all of the glossary entries
\usepackage{comment}
\usepackage[colorlinks]{hyperref}
\usepackage[nameinlink]{cleveref}

\hypersetup{
linkcolor=BrickRed
,citecolor=Green
,filecolor=Mulberry
,urlcolor=NavyBlue
,menucolor=BrickRed
,runcolor=Mulberry
,linkbordercolor=BrickRed
,citebordercolor=Green
,filebordercolor=Mulberry
,urlbordercolor=NavyBlue
,menubordercolor=BrickRed
,runbordercolor=Mulberry
}

\makeglossaries
\newacronym{ai}{AI}{artificial intelligence}
\newacronym{biv}{BiV}{biventricular}
\newacronym{bivp}{BiVP}{biventricular pacing}
\newacronym{crt}{CRT}{cardiac resynchronization therapy}
\newacronym{cs}{CS}{conduction system}
\newacronym{csp}{CSP}{conduction system pacing}
\newacronym{ivs}{IVS}{interventricular septum}
\newacronym{lbb}{LBB}{left bundle branch}
\newacronym{lbbb}{LBBB}{left bundle branch block}
\newacronym{lvsp}{LVSP}{left ventricular septal pacing}
\newacronym{hbp}{HBP}{His bundle pacing}
\newacronym{lbbap}{LBBAP}{left bundle branch area pacing}
\newacronym{lbbp}{LBBP}{left bundle branch pacing}
\newacronym{lv}{LV}{left ventricle}
\newacronym{nslbbp}{ns-LBBP}{non-selective left bundle branch pacing}
\newacronym{rv}{RV}{right ventricle}
\newacronym{slbbp}{s-LBBP}{selective left bundle branch pacing}
\newacronym{sd}{SD}{strength-duration}
\newacronym{vep}{VEP}{virtual electrode polarization}

\newcommand{\keypoints}{{\small
\begin{adjustwidth}{1cm}{1cm} 
    \centering{\textbf{Key Points}}
    \begin{itemize}
    \item Selectively capturing the \gls{lbb} requires direct physical contact between electrode tip and \gls{lbb}.
    \item Capturing the \gls{lbb} in its entirety may not always be feasible with lower selective stimulus strength.
    \item Stimulation with anodal polarity at the electrode tip requires a higher pulse strength for \gls{lbb} activation.
    \item The \gls{slbbp} capture threshold is influenced by the lead orientation relative to the \gls{lbb} bundles.
    \item Impedance trends observed clinically when advancing the lead through the \gls{ivs}
          are explained by the variability in tissue conductivity surrounding the helix in the \gls{lv} subendocardium,
          comprising myocardium, \gls{lbb} and fibrous tissue.
    \item A significant impedance drop, coupled with an increased \gls{slbbp} capture threshold is indicative of septal perforation.
    \item Lead dislodgement is indicated by an impedance drop without a significant change in the \gls{slbbp} capture threshold.
    \item Revealed mechanisms and quantitative consistence with clinical trends support model credibility, 
          suggesting the use of simulations as an effective approach for guiding the design of improved \gls{csp} leads.
    \end{itemize}
  \end{adjustwidth}}
  \vspace{2em}
}

\title{Computational Modeling of Selective Capture Mechanisms in Conduction System Pacing}
\author[1]{\small Mohammadreza Kariman} 
\author[1]{\small Matthias A.F. Gsell}
\author[4]{\small Edward J. Vigmond}
\author[3]{\small Aurel Neic} 
\author[1,2]{\small Christoph M. Augustin} 
\author[1,2]{\small Gernot Plank\thanks{%
\small Correspondence to  Gernot Plank: Gottfried Schatz Research Center: Division of Medical Physics and Biophysics, Medical University of Graz, Neue Stiftingtalstraße 6(MC2.H.)/III, 8010 Graz, Austria. Email: gernot.plank@medunigraz.at
}} 
\affil[1]{\small Gottfried Schatz Research Center: Medical Physics and Biophysics, Medical University of Graz, Graz, Austria} 
\affil[2]{\small BioTechMed-Graz, Graz, Austria} 
\affil[3]{\small NumeriCor GmbH, Graz, Austria}
\affil[4]{\small IHU Liryc, Electrophysiology and Heart Modeling Institute, Fondation Bordeaux Universit\'e, Liryc, Bordeaux, France}

\begin{document}
\onehalfspacing
\maketitle
\keypoints
\begin{abstract}
% Motivation & Background
\Gls{csp} is gaining clinical significance owing to its ability to restore a physiological activation sequence in the ventricles.
While \gls{hbp} producing the most physiological activation is preferable, 
due to implant complications the selective activation of the \gls{lbb} by \gls{lbbap} is considered an alternative,
offering both a simpler implant and a physiological activation sequence. 
However, the physical mechanisms facilitating selective activation of the \gls{lbb} remain poorly understood.
% Methods
We developed a structurally and biophysically detailed computer model of the \gls{ivs} and \gls{lbb}
to quantitatively elucidate the role of lead position, orientation and polarity in achieving optimal \gls{slbbp} thresholds,
using a geometrically detailed model of a clinically widely used \gls{csp} lead. 
% Results
A deep implant within the \gls{lv} sub-endocardium ensuring a direct contact between electrode and \gls{lbb} 
is key for effective \gls{slbbp}.
For low strength \gls{slbbp} is feasible, but capturing the \gls{lbb} in its entirety could only be achieved
using higher strengths that led to \gls{nslbbp}.
Switching the tip polarity to anodal was not beneficial, requiring higher strengths to activate the \gls{lbb}. 
Lead orientation relative to the \gls{lbb} bundles was found to influence the \gls{slbbp} capture threshold 
and the number of synchronously activating bundles.
% The model explains the impedance trends that are clinically observed 
% when advancing the tip through fibrous tissue and \gls{lbb},
The model explains the impedance trends that are clinically observed when advancing the tip through the \gls{ivs} into the \gls{lbb} region, 
as well as sudden impedance drops associated with implant complications such as septal perforation or lead dislodgement.
% Significance
% where the latter could be distinguished by a \gls{slbbp} capture threshold that remained nearly unchanged.
Quantitative consistence with clinically observed trends support model credibility, 
and indicate that simulation may offer an effective approach for guiding the design of improved \gls{csp} leads,
facilitating a selective and synchronous activation of the entire \gls{lbb}.

\end{abstract}

\paragraph{Keywords:} Conduction System Pacing; Left Bundle Branch Area Pacing; Computational Modelling; Cardiac Resynchronisation Therapy.
%

% ============================================= Introduction ================================================

\section{Introduction}\label{sec:introduction}

% Background & motivation
Ventricular dyssynchrony, resulting from \gls{cs} abnormalities or various cardiac diseases, 
significantly contributes to the loss of synchronous ventricular contraction, 
which is a key factor in the progression of heart failure~\cite{ghio_interventricular_2004}. 
\Gls{crt} with \gls{bivp} aiming at restoring a synchronized ventricular contraction that improves overall cardiac output
~\cite{auricchio_long-term_2002, moss_cardiac-resynchronization_2009}
%along with pharmacologic management 
has emerged as a primary treatment
for preventing progression or even achieving remodeling of heart failure in patients
with an ejection fraction of \SI{\leq 35}{\%} and an QRS duration longer than \SI{\geq 120}{\milli \second}~\cite{epstein_accahahrs_2008}.
%by restoring a synchronized ventricular contraction that improves overall cardiac output
%~\cite{auricchio_long-term_2002, moss_cardiac-resynchronization_2009}. 
However, despite sustained efforts to improve \gls{crt}~\cite{vernooy_strategies_2014, van_der_wall_improvement_2014},
the rate of non-responders that do not derive clinical or echocardiographic benefits remains high,
between \SI{30}{\%} and \SI{40}{\%}, and some may even worsen \cite{vernooy_strategies_2014, ruschitzka_cardiac-resynchronization_2013}.

% \gls{bivp} has been shown to improve cardiac function and survival in patients with an ejection fraction of $35\%$ or less and an intraventricular conduction delay of 120 ms or greater~\cite{epstein_accahahrs_2008}. Despite its effectiveness, the reported non-responder rate is as high as $30\%$ potentially due to the non-physiological fusion of stimulated wavefronts~\cite{varma_evaluation_2019, daubert_avoiding_2016}.  

While response to \gls{biv}-\gls{crt} is dependent on multiple factors 
such as the choice of pacing sites~\cite{parreira_non-invasive_2023, zweerink_hemodynamic_2019},
a fundamental limitation is the pacing-induced non-physiological activation sequence ~\cite{varma_evaluation_2019, daubert_avoiding_2016}
which differs markedly from the healthy intrinsic synchronous activation
and may, over time, contribute to the development of adverse effects~\cite{ruschitzka_cardiac-resynchronization_2013}.
Recently, this led to the emergence of \gls{csp} as a promising alternative to \gls{bivp}. 
%to deliver physiological ventricular pacing~\cite{deshmukh_permanent_2000}. 
\Gls{csp} involves stimulating the cardiac \gls{cs} at the His bundle or its primary bundle branches 
aiming at preserving a physiological intra-ventricular activation~\cite{deshmukh_permanent_2000}
and delivering \gls{crt} by synchronizing inter-ventricular activation~\cite{burri_ehra_2023, ali_left_2023}. 
%In the context of \gls{csp}, 
For \gls{csp} delivery by \gls{hbp} and \gls{lbbp} superior ventricular synchronization has been shown 
compared to \gls{bivp}~\cite{zhang_left_2019, vijayaraman_comparison_2023, sharma_permanent_2018}. 
However, achieving optimal \gls{csp} delivery that restores, as closely as possible, 
the intrinsic physiological ventricular activation sequence remains challenging
due to factors such as e.g.\ difficulties in lead positioning,
%sub-optimal electrical lead performance -- particularly the progressive rise in pacing thresholds --
%limitations in addressing infra-Hisian or more distal conduction disorders, 
and controlling capture to achieve optimal ventricular activation.
For restoring a most physiological activation sequence 
selective capture of the His or \gls{lbb} is considered optimal, 
but surrounding myocardium may also be excited, resulting in non-selective capture, 
which may result in different clinical outcomes depending on the patient's underlying cardiovascular condition~\cite{upadhyay_selective_2017, zhang_comparison_2018}.
Whether capture is \gls{slbbp} or \gls{nslbbp} is influenced by
pacing threshold parameters such as pulse strength and duration, 
lead implantation site and orientation,
differences in structural and electrophysiological properties~\cite{massing_anatomical_1976, stephenson_high_2017, padala_anatomy_2021}. 
and the variability of the cardiac \gls{cs}
across different anatomical regions~\cite{wu_evaluation_2021, burri_ehra_2023}.

% This variability encompasses both the structural and electrophysiological properties of the surrounding tissue, including conduction characteristics~\cite{massing_anatomical_1976, stephenson_high_2017, padala_anatomy_2021}. 

%These factors combined with the inherently higher complexity of the implant procedure
%hinder the widespread adoption of \gls{hbp} as a routine alternative to conventional \gls{biv}-\gls{crt} \cite{glikson_european_2025, sharma_permanent_2018, vijayaraman_his-bundle_2019}

% \gls{lbbp}, an approach within \gls{lbbap}, produces distinct physiological responses depending on whether the conduction system is activated selectively or non-selectively. 
% In \gls{slbbp}, the \gls{cs} is exclusively activated without direct myocardial capture, 
% whereas \gls{nslbbp} involves simultaneous activation of both the \gls{cs} and the adjacent myocardial tissue near the pacing site, 
% which may result in different clinical outcomes depending on the patient's underlying cardiovascular condition~\cite{upadhyay_selective_2017, zhang_comparison_2018}.

% These factors not only affect the transition between s-LBBP and ns-LBBP but also play a crucial role in determining patient-specific physiological responses to CSP. 

% \todo[inline]{Suggest to remove, we do not need this, maybe the citations is relevant?}

Furthermore, tissue-device interactions significantly impact the evolving application of this \gls{crt} strategy, a comprehensive evaluation of these variables is essential for optimizing therapeutic efficacy and improving clinical outcomes~\cite{sun_influence_2022}. 

% Consequently, the broader implementation of CSP in clinical practice is impeded by the absence of large-scale randomised trials validating its outcomes, coupled with persistent concerns about its long-term efficacy, especially in patients with heart failure, where the full scope of its benefits and potential risks has yet to be fully established~\cite{domenichini_conduction_2023}.

% In recent years, a number of computational studies have been conducted  to advance the understanding of CSP from different aspects \cite{vigmond_how_2021, strocchi_computational_2024, mirmaksudov_enhancing_2024}. 

Key to optimizing \gls{csp} is a better mechanistic understanding of the physics 
governing device-tissue interaction
by which the \gls{cs} and tissue adjacent to a \gls{csp} lead is captured. 
In this study, we use a biophysically detailed computational model of the \gls{ivs}
and the \gls{lbb} along with a geometrically accurate representation of a clinically widely used \gls{csp} lead.
By simulating the spatio-temporal evolution of the extracellular potential field 
imposed by the lead during stimulation and the evoked polarization within \gls{lbb} and \gls{ivs} myocardium,
capture mechanisms and differences in capture thresholds are investigated.
The role of most relevant parameters like lead position relative to the \gls{lbb}, 
orientation, polarity and the applied simulation protocol is quantitatively analyzed.
By simulating lead advancement, model predictions are validated against clinical measurements 
and observable trends, including during known implant complications 
such as septal perforation or lead dislodgement.
The achieved quantitative accuracy and the consistent prediction of clinically observed trends suggest
that our approach is suitable for aiding in \gls{csp} lead design for robustly achieving \gls{slbbp}. 

%This study aimed to offer a new perspective on the CSP mechanism, specifically LBBAP, 
%by investigating tissue-device interactions through computational modelling. 
%The focus was placed on the detailed examination of selective and non-selective LBBP capture mechanisms, as well as potential complications at the tissue-device interface. 
%By offering a novel insight into these mechanisms from a fresh outlook, this work contributes to a deeper understanding of this evolving CRT approach, supporting its future clinical refinement. 

% ================================================ Method ===================================================

\section{Methods}\label{sec:method}

\subsection{Computational Model of Tissue and Device}
A computational model of the \gls{ivs} was generated,
with dimensions of \SI{5}{\centi \meter} $\times$ \SI{5}{\centi \meter} $\times$ \SI{1}{\centi \meter}.
Septal fiber architecture was incorporated by employing a rule-based approach \cite{bayer_novel_2012},
imposing a fiber rotation from $+60^\circ$ to $-60^\circ$ 
on \gls{lv} and \gls{rv} endocardial surfaces respectively, 
in accordance histological studies \cite{streeter_fiber_1969}. 
The non-branching segment of the \gls{lbb} was modeled as $20$ parallel bundles, 
each with a diameter of \SI{200}{\micro \meter}, 
collectively forming a \SI{\approx 10}{\milli\meter}-wide band-like structure  
consistent with reported dimensions \cite{massing_anatomical_1976}. 
These cables were positioned within a subendocardial layer of \SI{100}{\micro \meter} width 
and \SI{200}{\micro \meter} below the \gls{lv} endocardial surface \cite{cabrera_tracking_2020}.
The cables forming the \gls{lbb} were electrically not connected amongst each other, 
and surrounded by a thin layer of fibrous tissue 
which electrically insulated the \gls{lbb} fibers from the adjacent myocardial tissue of the \gls{ivs}.
Sufficiently large blood pools of \SI{5}{\centi \meter} depth
were attached to both \gls{lv} and \gls{rv} endocardia of the \gls{ivs}
to ascertain for all lead deployment configurations
that the blood pools covered the entire lead,
maintaining a minimum distance to the blood pool boundary of \SI{>1}{\centi \meter}.
Domains were conformally discretized at average spatial resolutions of 
\SI{250}{\micro \meter}, \SI{25}{\micro \meter} and \SI{500}{\micro \meter} 
for \gls{ivs}, \gls{lbb} cables and blood pools, respectively.

%The entire structure was subsequently embedded within a larger volume representing the blood pool. 
A model of a pacing lead was developed based on a commercially available device 
(Medtronic SelectSecure 3830) that is widely used for \gls{csp} clinically. 
Lead geometry and physical properties were based on the specification data sheet \cite{medtronic3830}.
The lead geometry was discretized at an average mesh resolution of \SI{100}{\micro \meter}, 
with a refined resolution of \SI{50}{\micro \meter} at the helical tip,
and integrated into the \gls{ivs} in various deployment configurations.
All processing steps involving device positioning, 
volumetric remeshing including adaptive conformal mesh resolution adjustments 
between device and model domains, were automated using \emph{Studio} 
(\href{https://numericor.at}{NumeriCor GmbH}, Graz, Austria).

\subsection{Simulating Electrophysiology and Electrical Stimulation}
Electrophysiology in \gls{ivs} and \gls{lbb} as well as associated current flow in the surrounding medium due to both 
electric potential field imposed by the device as well as electrical activity of the tissue 
was represented by the  mechanistically most comprehensive bidomain model~\cite{vigmond08:_solvers}. 
Cellular dynamics in the \gls{ivs} and \gls{lbb} were represented by models of ventricular and Purkinje myocytes proposed by Ten Tusscher \cite{tusscher_cell_2006} and Stewart \cite{stewart_mathematical_2009}, respectively.

%This combined approach allowed for an in-depth exploration of the electrical activity within both myocardial and conduction system tissues \cite{tusscher_cell_2006, stewart_mathematical_2009}.

% \todo[inline]{Moh, please document all settings used in tunecv here as a comment.}

Conductivities within the \gls{ivs} were calibrated based on the mesh resolution 
to attain conduction velocities of 
\SI{0.6}{\meter / \second}, \SI{0.4}{\meter / \second} and \SI{0.2}{\meter / \second} 
along the fiber, sheet and sheet normal direction, respectively,
yielding intra- and extracellular conductivities of 
$\sigma_{\mathrm{i,f}}$=\(0.157 \, \mathrm{S/m}\) and $\sigma_{\mathrm{e,f}}$=\(0.62 \, \mathrm{S/m}\),
$\sigma_{\mathrm{i,s}}$=\(0.076 \, \mathrm{S/m}\) and $\sigma_{\mathrm{e,s}}$=\(0.24 \, \mathrm{S/m}\),
$\sigma_{\mathrm{i,n}}$=\(0.02 \, \mathrm{S/m}\) and $\sigma_{\mathrm{e,n}}$=\(0.24 \, \mathrm{S/m}\), 
reflecting physiologically relevant anisotropy \cite{https://doi.org/10.1113/jphysiol.1976.sp011283}. 
For the \gls{lbb}, an isotropic conduction velocity of \(2.0 \, \mathrm{m/s}\) was considered in all directions,
yielding an intra- and extracellular conductivity of 
$\sigma_{\mathrm{i}}$=\(0.157 \, \mathrm{S/m}\) and $\sigma_{\mathrm{e}}$=\(0.62 \, \mathrm{S/m}\).
% with a four-fold increase in peak sodium conductance and a surface-to-volume ratio of \SI{1000}{\centi \meter}$^{-1}$.

fibrous tissue and blood pool were assigned conductivities of 
\SI{0.04}{S/\meter} and \SI{0.7}{S/\meter}, respectively~\cite{miklavvcivc2006electric}. 
Tip and ring electrode, assumed to be made of platinum, were assigned a conductivity of 
\( 9 \times 10^{6} \, \mathrm{S/m} \) \cite{flynn_measurements_1967},
whereas for the lead body, assumed to be composed of polyurethane, a conductivity of 
\( 10^{-11} \, \mathrm{S/m} \) was assigned~\cite{jafarzadeh_review_2023}.

All simulations were conducted using the finite element framework, \emph{Cardiac Arrhythmia Research Package (CARPentry)}\cite{vigmond08:_solvers}, which is built upon extensions of the \emph{openCARP} framework \cite{plank_opencarp_2021} (\url{http://www.opencarp.org}).
The bidomain equations were solved using a semi-implicit Crank-Nicolson scheme
with a time step of \SI{25}{\micro \second}. 
All simulations were performed on the Vienna Scientific Cluster (VSC-4) utilizing 8 nodes, each equipped with 48 cores (\url{https://vsc.ac.at}).

\subsection{Role of lead deployment}

Lead deployment was systematically varied to investigate the effects of both lead proximity and orientation relative to the anatomical embedding of the \gls{lbb} within the \gls{ivs}.

Both factors are  major determinants of field-tissue interaction driving the change in transmembrane voltage, $V_{\rm m}$.
During delivery of a pacing stimulus, 
the deployed lead establishes an electric field, $\phi^p_{\rm e}$, given by
\begin{equation}
   \nabla \cdot \boldsymbol{\sigma} \nabla \phi^p_{\rm e} = 0, \label{eq:_bidm_ell}
\end{equation}
with Dirichlet boundary conditions at the helical tip, $\phi_{\rm t}$, and the ring electrode, $\phi_{\rm r}$, 
\begin{equation}
    \phi^p_{\rm e}|_{\rm {tip}} = \phi_{\rm t}(t)  \cup  \phi^p_{\rm e}|_{\rm {ring}} = \phi_{\rm r}(t). \label{eq:_bidm_bcond}
\end{equation}
Assuming that the potential field during stimulation is primarily dictated by the 
applied pacing field, $\phi_{\rm e}(\mathbf{x}) \approx \phi^p_{\rm e}(\mathbf{x})$,
the change in $V_{\rm m}$ over time is given then by
\begin{equation}
    C_{\rm m} \frac{\partial V_{\rm m}}{\partial t} = \underbrace{\nabla \cdot \boldsymbol{\sigma}_{\rm i} \nabla \phi^p_{\rm e}}_{S} -I_{\rm{ion}} + \nabla \cdot \boldsymbol{\sigma}_{\rm i} \nabla V_{\rm m}, \label{eq:_bidm_parab}
\end{equation}
where $S$ is referred to as activating function \cite{sobie_1997_generalized,rattay_analysis_1986}. 
The strength of the field contribution of $\phi^p_{\rm e}$ to $S$ and, thus, $\Delta V_{\rm m}$, 
depends on the magnitude of the electric field weighted with the tissue heterogeneity 
$(\nabla \cdot \boldsymbol{\sigma}_{\rm i}) \cdot \nabla \phi^p_{\rm e}$, 
and the heterogeneity of the electric field weighted with the tissue conductivity, $\boldsymbol{\sigma}_{\rm i} : \nabla \cdot \nabla \phi^p_{\rm e}$.

\subsubsection{Lead position in \glsentryshort{lbbap}}   % use \glsentrytext for the full term and \glsentryshort for the abbreviation

\begin{figure}[H]
    \centering
    \includegraphics[width=1.0\linewidth]{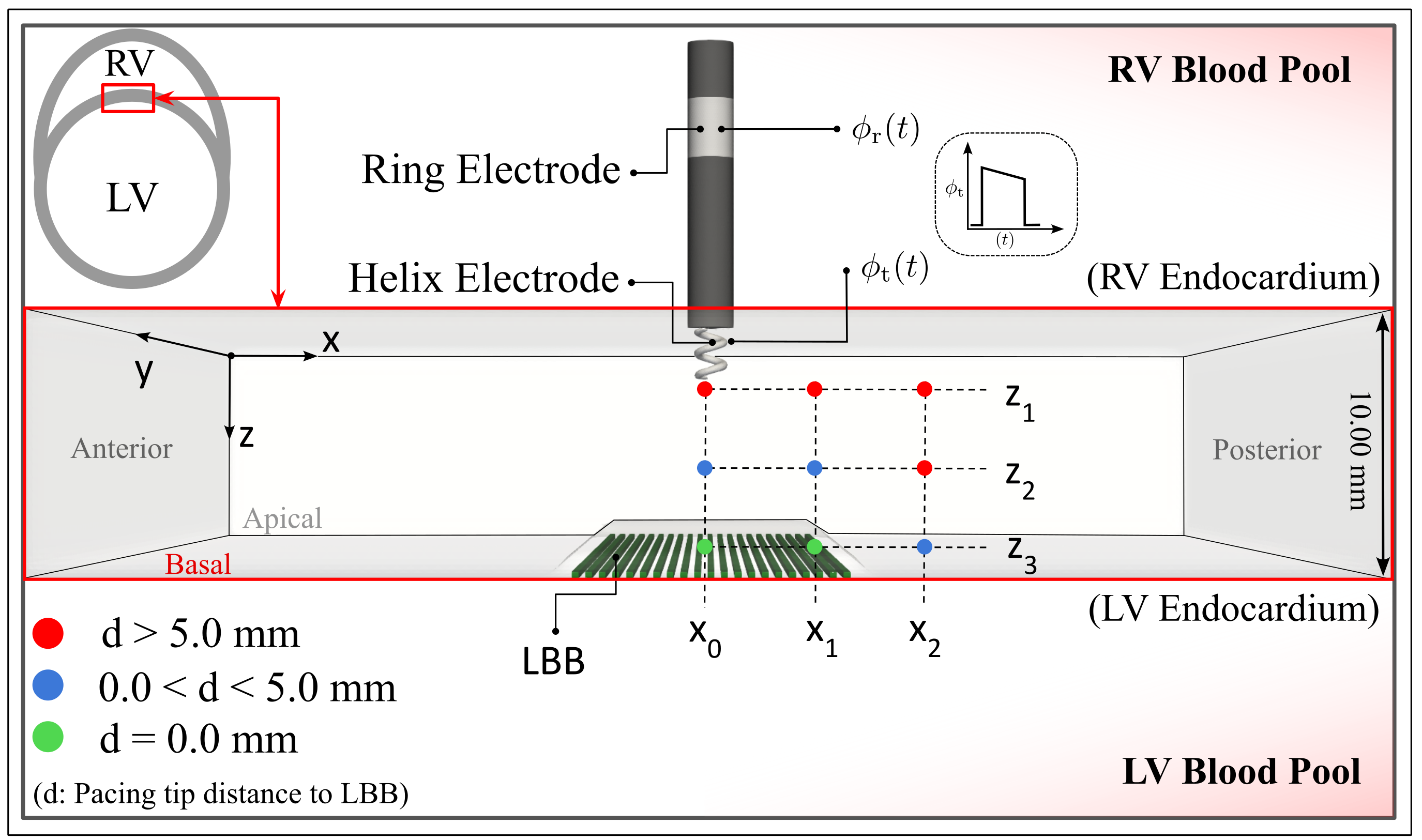}
    \caption{Simulation setup with \gls{csp} lead deployment in the \gls{ivs} 
    through an \gls{rv} access. 
    Model of an \gls{ivs} wedge with embedded \gls{lbb} and attached blood pools.
    Antero-posterior and transmural direction represented by $x$- and $z$-axis 
    encode the lead tip position within the \gls{ivs} at a given apico-basal height denoted by $y$.   
    Point colors encode the Euclidean distance between the tip and and \gls{lbb} 
    of all nine deployment positions.
    }
    \label{fig:setup}
\end{figure}

The electrophysiological response of the \gls{ivs} myocardium and \gls{lbb} to electrical stimulation via \gls{lbbap} was investigated by orthogonally advancing a pacing lead through \gls{rv} access into the septal endocardium (Figure~\ref{fig:setup}).

% Three deployments axis relative to the \gls{lbb} were tested, 
% starting at the midpoint of the \gls{lbb}, with two additional deployment axes positioned progressively farther posteriorly.
Three deployment points along an axis relative to the \gls{lbb} were tested, starting at the midpoint of the \gls{lbb}, with two additional points positioned progressively farther posteriorly.
These deployment axes are represented by the coordinates \(x_0\), \(x_1\), and \(x_2\) in Figure \ref{fig:setup}. 
For each deployment axis, the lead was advanced to three different implantation depths, $z$:
\begin{itemize}
    \item Depth $z_1$: A superficial implantation corresponding to the helix penetration depth of \SI{1.8}{\milli \meter},
    \item Depth $z_2$: A mid-septal implantation at \SI{5}{\milli \meter}, positioning the  tip approximately mid-septally,
    \item Depth $z_3$: A deep implantation at \SI{10}{\milli \meter} where the tip 
    is in direct contact with the \gls{lbb},
    \end{itemize}
yielding nine different deployment locations with varying lead-to-\gls{lbb} distance in transmural and antero-posterior directions.

\subsubsection{Lead Orientation in \glsentryshort{lbbap}}

While for \gls{lbbap} a deployment perpendicular to the \gls{ivs} is recommended 
for reasons of efficient implantation achieving a deep and stable anchoring of the tip within the septum \cite{de_pooter_guide_2022, shroff_comparison_2024},
in practice, this cannot be precisely controlled or may not be achievable 
due to factors such as \gls{ivs} anatomy or restrictions imposed by the tricuspid valve. 
The influence of a non-orthogonal deployment was therefore further assessed for the tip position $(x_0, z_3)$ 
by inclining the lead axis from an orthogonal deployment by 20$^\circ$ 
along and transverse to the \gls{lbb} axis (Figure \ref{fig:advance_stories}: A \& B). 

Both inclinations rotate the electric field, $\mathbf{E}=-\nabla \phi^p_{\rm e}$,
to enhance the field oriented parallel or transverse to the \gls{lbb} compared to a purely orthogonal deployment.   

% \begin{itemize}
%     \item Inclination Along the LBB Axis: The lead was inclined $20^\circ$ towards the long axis of the LBB, positioning the ring electrode closer to its proximal region.
%     \item Inclination Towards the LBB Transverse Axis: The lead was inclined $20^\circ$ towards the LBB's transverse axis, causing the ring electrode to rotate posteriorly.
%\end{itemize}
%These configurations were designed to evaluate the impact of azimuthal orientation on optimizing LBB electrical capture in computational simulations.

\subsubsection{Septal Perforation}
During implantation, the lead is gradually advanced deeper into the \gls{ivs} to position the helical tip as close as possible to, or within, the \gls{lbb} band. If advanced further, the tip electrode may perforate the \gls{lv} endocardium, resulting in contact with the \gls{lv} blood pool.
We investigate the effect of septal perforation on \gls{lbb} capture, assuming full penetration such that the helical tip protrudes entirely into the \gls{lv} blood pool (Fig. \ref{fig:advance_stories}). 

\subsubsection{Lead Dislodgement}
Lead dislodgement is a potential complication during \gls{lbbap}, caused by a failure of the drilling mechanism that impedes further advancement of the lead into the \gls{ivs} tissue. This results in a tunnel-like structure within the \gls{ivs}, hindering firm anchoring of the lead~\cite{burri_ehra_2023}.

\begin{figure}[H]
    \centering
    \includegraphics[width=0.5\linewidth]{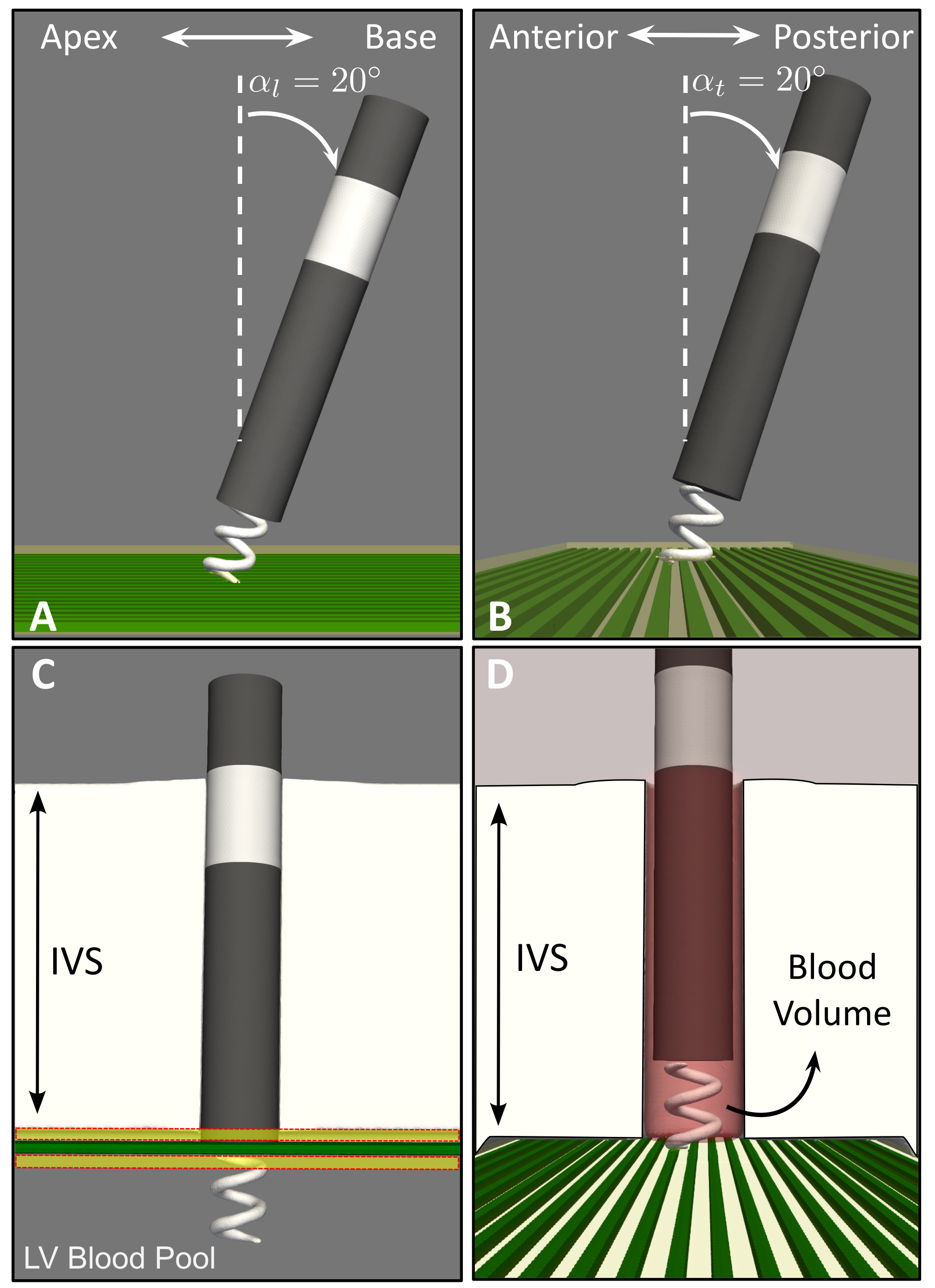}
    \caption{Variability and complications with lead deployment with \gls{lbbap}. 
    The lead axis is advanced along a direction non-orthogonal to the \gls{lbb} 
    under an inclination angle of 20$^\circ$ along the \gls{lbb} axis (A) or transverse to it (B). 
    (C) Septal perforation was modeled by advancing the lead through the \gls{ivs} 
    until the helical tip penetrates the insulating fibrous tissue (yellow layer)
    and protrudes into the \gls{lv} blood pool. 
    (D) Lead dislodgement was modeled by creating a blood-filled lumen around the helical tip.   
    }
    \label{fig:advance_stories}
\end{figure}

The effect of lead dislodgement was modeled for orthogonal deployment at location $(x_0,z_3)$ 
by adding a blood-filled lumen within the \gls{ivs}
surrounding the helical tip (Fig.~\ref{fig:advance_stories}D). 

%This approach aimed to improve the understanding of the response dynamics during lead dislodgement and its subsequent impact on LBBP performance.

%\todo[inline]{I'd suggest to combine Fig 1. and 2., all of this belongs together.}

\subsection{Validation of model predictions}
\label{sec:validation}
% rather results than methods
Clinical measurements routinely performed during \gls{csp} implant 
include capture thresholds, indicators of selective versus non-selective capture and monitoring of impedances.
As these measurements are highly variable between patients, 
and directly measured data were not available, 
we refrained from calibration of model parameters to improve the goodness of fit. 
Rather, tissue and lead parameters based on values reported in the literature were used
to compute capture thresholds, lead impedances and type of capture. 
Predictions were validated then qualitatively by assessing 
whether computed values fall within the range of clinical observations, 
and whether the predicted trends are consistent with clinical experience.

\subsubsection{\glsentryshort{sd} Curves}\label{sec:SD_Curves}
For each deployment configuration stimulus \gls{sd} curves were computed 
to determine the excitation thresholds, $\varepsilon_{\rm{ivs}}$ and $\varepsilon_{\rm{lbb}}$, 
for \gls{ivs} myocardium and \gls{lbb}, respectively.

The delivery of monophasic rectangular pulses of variable strength in voltage, $\Delta \phi_{\mathrm{stim}}$,
imposed between helical tip and ring, $\Delta \phi_{\mathrm{stim}} = \phi_{\rm t} - \phi_{\rm r}$, 
polarity and duration $\tau$ was simulated.
Excitation thresholds, $\varepsilon_{\rm{ivs}}$ and $\varepsilon_{\rm{lbb}}$, were estimated 
with an accuracy of \SI{\pm 0.05}{\volt} using a bisection method  
for finding the root of $\Delta \phi_{\mathrm{stim}}(\tau) - \varepsilon = 0$. 
For each duration $\tau$, 
sampled in the range from $\tau_{\rm {min}}=$ \SI{0.5}{\milli \second} to $\tau_{\rm {max}}=$ \SI{4}{\milli \second} in increments of \SI{0.5}{\milli \second},
variable stimulus strength in the range from $\Delta \phi_{\rm {stim,min}}=$ \SI{0.1}{\volt} to $\Delta \phi_{\rm {stim,max}}=$ \SI{15}{\volt} were tested,
to determine the corresponding thresholds $\varepsilon_{\rm{ivs}}(\tau)$ and $\varepsilon_{\rm{lbb}}(\tau)$.

For each test $\Delta \phi_{\mathrm{stim}}(\tau)$ the response to simulation was classified into four categories:
%A python algorithm was developed to manage the simulations required for constructing accurate strength–duration (S-D) curves in LBBP. 
%The approach was based on the bisection method, a numerical technique commonly used to find roots of functions. 
%In this context, it was applied to determine the precise stimulus parameters—specifically amplitude and pulse duration—necessary to achieve distinct activation outcomes. This method aimed to accurately identify the boundaries where transitions occur between different physiological responses, particularly between selective and non-selective capture.

%By iteratively refining stimulus parameters within a predefined range, threshold identification and optimal stimulation settings were ensured and the outcomes of these simulations were classified into four categories:

\begin{enumerate}
    \item \textbf{Loss of Capture}: The stimulus $\Delta \phi_{\mathrm{stim}}(\tau)$ was insufficient, 
    neither \gls{ivs} nor \gls{lbb} were activated, 
    $\Delta \phi_{\mathrm{stim}}(\tau) < \varepsilon_{\rm {ivs}}$ and $\Delta \phi_{\mathrm{stim}}(\tau) < \varepsilon_{\rm {lbb}}$,
    \item \textbf{Selective Capture}: Only the \gls{lbb} was activated, but not the \gls{ivs}, 
    $\Delta \phi_{\mathrm{stim}}(\tau) < \varepsilon_{\rm {ivs}}$ and $\Delta \phi_{\mathrm{stim}}(\tau) >= \varepsilon_{\rm {lbb}}$,
    \item \textbf{Non-Selective Capture}: Both \gls{lbb} and \gls{ivs} were activated,
    $\Delta \phi_{\mathrm{stim}}(\tau) >= \varepsilon_{\rm {ivs}}$ and $\Delta \phi_{\mathrm{stim}}(\tau) >= \varepsilon_{\rm {lbb}}$,
    \item \textbf{Myocardial Capture}: The \gls{ivs} was activated without engaging the \gls{lbb},
    $\Delta \phi_{\mathrm{stim}}(\tau) >= \varepsilon_{\rm {ivs}}$ and $\Delta \phi_{\mathrm{stim}}(\tau) < \varepsilon_{\rm {lbb}}$,
\end{enumerate}
where the decision boundaries $\varepsilon(\tau)$ separated the domains with and without capture.

%The simulation process began by selecting the midpoint of a predefined stimulus amplitude range, using the shortest pulse duration (0.5 ms in this study). Subsequent simulations iteratively refined the stimulation threshold through a bisection search algorithm, which halved the test range based on each outcome until the threshold was accurately determined. Once the threshold for a given pulse duration was identified, the search range for subsequent durations was adaptively narrowed. Previously computed results were leveraged to enhance computational efficiency, reducing the number of iterations required for later calculations. This methodology enabled the precise mapping of the S-D curve, providing a robust framework for optimizing stimulation parameters in conduction system pacing. 

%As an example, to identify the lowest s-LBBP threshold, the process was halted for a given pacing duration when the difference between the highest tested amplitude that resulted in loss of capture and the lowest tested amplitude yielding selective capture fell below the predefined tolerance

To account for the unknown variability in the tissue compositions around the electrode tip 
additional \gls{sd} curves were computed to investigate the dependency of capture thresholds with \gls{lbbp}.
Specifically, the following scenarios were considered:
\begin{enumerate}
    \item The pacing lead was implanted in the \gls{ivs} and positioned near the \gls{lv} subendocardial layer, 
    where most of the tip surface area (70\%) was covered by fibrous tissue surrounding the \gls{lbb}, 
    and the remaining surface area was covered to 25\% by \gls{ivs} myocardium (25\%) and to 5\% by the \gls{lbb}.

    \item The helix, with no connection to the \gls{lbb}, maintained contact with myocardium (25\%) and fibrous tissue (75\%).

    \item The helix, with no connection to myocardium, retained its connection with the \gls{lbb} (5\%) and fibrous tissue (95\%).

\end{enumerate}

\subsubsection{Device Impedance}\label{sec:device_impedance}

Unipolar device impedance between helical tip and ring for each deployment configuration
computed by solving the elliptic portion of the bidomain equation \eqref{eq:_bidm_ell} for $\phi_{\mathrm{e}}(\mathbf{x})$,  
and determining the current density, $\mathbf{J}$, then by 
\begin{equation}\label{eq:J}
\mathbf{J}(\mathbf{x},t) = -\boldsymbol{\sigma}_\mathrm{e}(\mathbf{x}) \nabla \phi_\mathrm{e}(\mathbf{x,t}).
\end{equation}

The total current flowing over the tip, $I_{\mathrm{tip}}$, was calculated 
by integrating $\mathbf{J}$ over the tip surface, $\Gamma$, of the electrode 
over the duration $\tau$ of the stimulus pulse, and averaging it over $\tau$, 
yielding
\begin{equation}\label{eq:I}
I_{\mathrm{tip}} = \frac{1}{T} \int_{t=0}^{\tau} \int_{\Gamma_{\rm {tip}}} \mathbf{J}(\mathbf{x},t) \, d\boldsymbol{\Gamma} \, dt.
\end{equation}
The impedance $Z$ is calculated then by 
\begin{equation}\label{eq:Z}
Z = \frac{V}{I_{\mathrm{tip}}}
\end{equation}
where the voltage is determined as the difference between the potentials prescribed at tip and ring electrode, 
$V = \phi_{t} - \phi_{r}$.

\subsubsection{Implant Impedance Monitoring}\label{sec:implant_impedance}

Device impedance $Z$ during implant is variable as the lead is advanced deeper into the \gls{ivs}.
This variation of $Z$ with depth $z$ is continuously monitored 
as it is used for determining the optimal lead position 
and for detecting potential procedure-related complications~\cite{burri_ehra_2023}.
To simulate the deployment process, unipolar tip impedances $Z$ were computed 
at various depths $z$, and the course of $Z(z)$ compared to trends recorded during implant.

This course of $Z(z)$ aids in positioning the tip close to the \gls{lv} endocardium 
where the \gls{lbb} resides \cite{vernooy_implant_2023},
and in detecting conditions such as the perforation of the \gls{lv} endocardium or lead dislodgement.
These are detected by a drop in $Z$
as the tip will be covered by a highly conductive blood volume in both scenarios.

\subsubsection{Impact of tissue composition upon impedance}
Since unipolar tip impedance $Z$ is used for finding an optimal placement location,
and, potentially, for increasing the safety margin for \gls{slbbp},
better understanding the dependency of $Z$ upon tissue composition may be important 
for guiding both lead design as well as the implant procedure, 
but also to infer information on the tissue composition around the lead.
While the dimension of anatomical and structural entities representing \gls{ivs}, \gls{lbb}
and fibrous sheath are based on data reported in the experimental literature \cite{james_fine_1971, elizari_normal_2017, titus_normal_1973, uhley_visualization_1959, demoulin_histopathological_1972, garciabustos_quantitative_2017, vijayaraman_prospective_2019},
the variability in terms of geometry and electrical conductivities is large,
which introduces uncertainties in the estimates of $Z$ in our model.

We evaluate therefore the relative importance of the volumetric distribution 
and the associated heterogeneous conductive milieu the tip is exposed to 
when being advanced through the \gls{lbb},  
by systematically varying the tissue composition $\boldsymbol{\sigma}(\mathbf{x})$,
% (Fig.~\ref{fig:impedance}C), 
comprising fibrous tissue, \gls{ivs} myocardium, \gls{lbb} and blood, respectively.
The relationship $Z=f(\boldsymbol{\sigma}(\mathbf{x}))$ was quantified using Pearson’s correlation.
Finally, \gls{sd} curves were also computed for various tissue configurations 
to demonstrate that trends observed with the baseline configuration are preserved.

% =============================================== Results ===================================================

\section{Results}\label{sec:results}

\subsection{Impact of Lead Positioning on \glsentryshort{lbbap} Capture Mechanisms}
 
For $3 \times 3$ implantation sites differing in transmural depth $(z_1,z_2,z_3)$ 
and anterior-posterior displacement relative to the center of the \gls{lbb} $(x_0,x_1,x_2)$ (Fig.~\ref{fig:setup}) \gls{sd} curves were computed. 

%The major observations were classified as follows:
% this text is too verbose for its message, I simplified this.
% It was observed that, at a shallow lead deployment depth of 1.8 mm—corresponding to the height of the pacing electrode—neither s-LBBP nor ns-LBBP could be achieved within the defined stimulation threshold range of 0.1 to 5.0 V and 0.5 to 4.0 ms. 
% In all these cases, myocardial capture was the only observed activation mode and physiological response, 
% with a minimum threshold of 0.25 V at 0.5 ms.

For stimuli ranging from \SI{0.1}{} to \SI{5.}{\volt} for durations of \SI{0.5}{} to \SI{4}{\milli \second}
none of the shallow deployment locations at $z_1$ were able to activate the \gls{lbb}, 
only myocardial capture was observed, with a minimum threshold of $V=$ \SI{0.25}{\milli \volt} at $\tau=$ \SI{0.5}{\milli \second}. Similarly, for the same range of stimuli at the deeper mid-septal deployment locations $z_2$ (Fig.~\ref{fig:setup}), also no \gls{lbb} capture was observed,
% this drop is so small, would suggest to ignore
only myocardial capture was achieved, 
with a marginal drop in threshold voltage as compared to a deployment at $z_1$ (Fig.~\ref{fig:SD_Curves_1}A).
However, when testing substantially higher stimulus strengths 
in the range of up to \SI{15}{\milli \volt} \gls{lbb} could be captured non-selectively
for tip locations $x_1$ and $x_2$ above the \gls{lbb} (Fig.~\ref{fig:SD_Curves_1}B).

%The effect of lead implantation depth in LBBAP was further investigated by extending the lead depth from 1.8 mm to 5.0 mm (points X0Z1, X1Z1, and X2Z1 in Figure \ref{fig:setup}). 
% Similar to superficial implantation, no LBB activation was observed across these points within the threshold test range up to 5.0 V. However, the minimum capture threshold for myocardial activation slightly increased from 0.25 V to 0.27 V at 0.5 ms (Fig.~\ref{fig:SD_Curves_1}:A). To determine if LBB activation could be achieved at this depth, which corresponded to approximately halfway through the IVS, the stimulation threshold range was increased from 0.1–5.0 V to 5.0–15.0 V, while maintaining pulse durations between 0.5 and 4.0 ms.
%At positions closer to the LBB (X0Z1 and X1Z1 in Figure \ref{fig:setup}), 
%where the lead was positioned more medially along the antero-posterior axis, 
%ns-LBBP was successfully achieved by increasing the stimulation threshold. 
%The LBB was activated non-selectively with a minimum threshold of 13.0 V at 0.5 ms when the pacing electrode tip was positioned within a Euclidean distance of less than 5 mm from the LBB (Figure \ref{fig:SD_Curves_1}:B).

Further advancing the lead to depth $z_3$ positioned the tip within the subendocardial layer where the \gls{lbb} is embedded. At the two positions along the antero-posterior axis where the tip made direct contact with the \gls{lbb} ($x_0$ and $x_1$ in Fig.~\ref{fig:setup}), \gls{slbbp} was observed.
However, when the tip was positioned at $x_2$, still within the subendocardial layer but spatially offset from the \gls{lbb} band, \gls{slbbp} could not be achieved. At positions where \gls{slbbp} was achieved, increasing the stimulation strength by a small margin of \SI{0.06}{\volt} for the same pulse duration resulted in a transition to \gls{nslbbp} (Fig.~\ref{fig:SD_Curves_1}C).

% \todo[inline]{I am thinking of adding a figure showing the space around the lbb 
where selective pacing could be achieved, and color coding every where the delta in Volts 
to the ns-lbbp threshold, may require a lot of computing though, but should show us the best spot
where the most robust s-lbbp could be achieved.

\subsection{Impact of Polarity - Anodal Pacing}
In conventional \gls{csp} -- whether using \gls{hbp} or \gls{lbbap} -- 
the helical tip electrode is configured as the cathode.

\begin{figure}[H]
    \centering
    \includegraphics[width=1.0\linewidth]{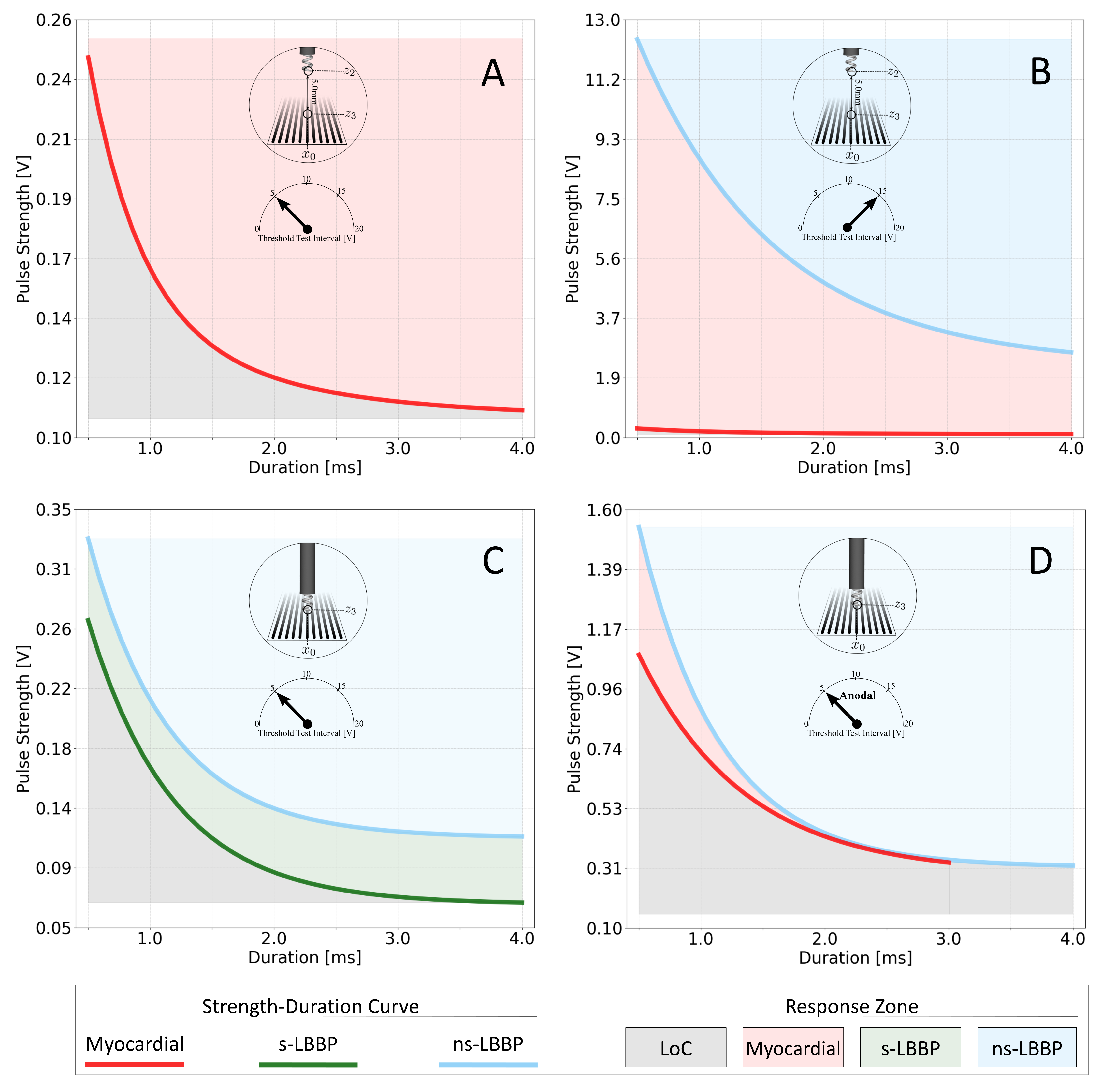}
    \caption{ \Gls{sd} curves for varying implantation depths $z$ and stimulation polarity.
    All implant sites were centered above the central axis of the \gls{lbb} at $x_0$. 
    Shown are \gls{sd} curves for:
    (A) depth $z_2$, resulting in myocardial-only capture;
    (B) depth $z_2$ with higher stimulus strengths up to \SI{15}{\volt}, resulting in \gls{nslbbp};
    % if the following is important, add the number to the figure!
    % at 12.8 V and 0.5 ms pulse width. 
    (C) depth $z_3$ where the tip is in contact with the \gls{lbb}, resulting in both
    \gls{slbbp} and \gls{nslbbp};
    % If you want these specific numbers, add them to the figure! All this information is in your graph, 
    % no need to repeat this in the caption!
    %(C) Lead implanted at a transmural depth of 9.0 mm with the helix making optimal contact with the LBB, resulting in s-LBBP within the tested threshold up to 5.0 V. s-LBBP achieved at 0.27 V @ 0.5 ms, 
    %while ns-LBBP occurred at 0.33 V @ 0.5 ms. 
    (D) altered tip polarity from cathode to anode, 
    resulting in an increase in capture threshold for both myocardium and \gls{lbb}
    and the loss of \gls{slbbp}.}
    %myocardial capture at 1.15 V @ 0.5 ms, and ns-LBBP at 1.6 V @ 0.5 ms.}
    \label{fig:SD_Curves_1}
\end{figure}

As anodal pacing captures tissue by a different mechanism, the \gls{sd} curves were computed for anodal pacing over the same stimulus strengths as for cathodal pacing to evaluate its relative efficacy for \gls{lbbap}. \Gls{sd} curves reveal that both the myocardium and the \gls{lbb} feature a higher capture threshold (compare Figs.~\ref{fig:SD_Curves_1}C and D).
Upon switching from cathodal to anodal pacing
the capture thresholds for a \SI{0.5}{\milli \second} stimulus increased, 
for the myocardium from \SI{0.25}{\volt} to \SI{1.5}{\volt} 
and for the \gls{lbb} from \SI{0.33}{\volt} to \SI{1.6}{\volt} (Fig.~\ref{fig:SD_Curves_1}:D).
As both \gls{sd} curves were shifted upwards towards higher strengths 
only myocardial and \gls{nslbbp} capture was observed, but not \gls{slbbp} capture. 
Mechanistically, capture of the \gls{ivs} myocardium and wave front initiation differed between the polarities. While cathodal stimulation resulted in tissue activation exclusively by the tip electrode, anodal stimulation -- requiring higher strengths -- led to wavefront initiation from both the tip and ring electrodes.

% The \gls{sd} curve resulting from simulations with the helix as the anode, conducted within the same threshold test range as cathodal excitation, revealed that both the myocardium and the LBB require a stronger pulse strength to be activated under anodal excitation.
% The increased pulse strength required for anodal excitation led to a loss of selective capture, 
% resulting in two physiological responses during the threshold tests: 
% myocardial-only and \gls{nslbbp} responses. 
% The capture threshold for the myocardium increased from 0.25 V to 1.15 V upon switching from cathodal to anodal excitation, while the threshold for non-selective activation rose from 0.33 V to 1.6 V, with both thresholds performed at 0.5 ms (Figure \ref{fig:SD_Curves_1}:D). 
% To further emphasize the differences, ns-LBBP resulting from cathodal excitation led to activation of the IVS exclusively by the helix electrode. In contrast, ns-LBBP resulting from anodal excitation caused activation of the IVS by both the helix and the ring electrode. This shift in activation patterns was driven by the significant increase in the capture threshold observed under anodal stimulation, despite both excitation types involving deep lead placement.

\subsection{\glsentryshort{lbb} Capture Mechanisms}
\begin{figure}[H]
    \centering
    \includegraphics[width=1\linewidth]{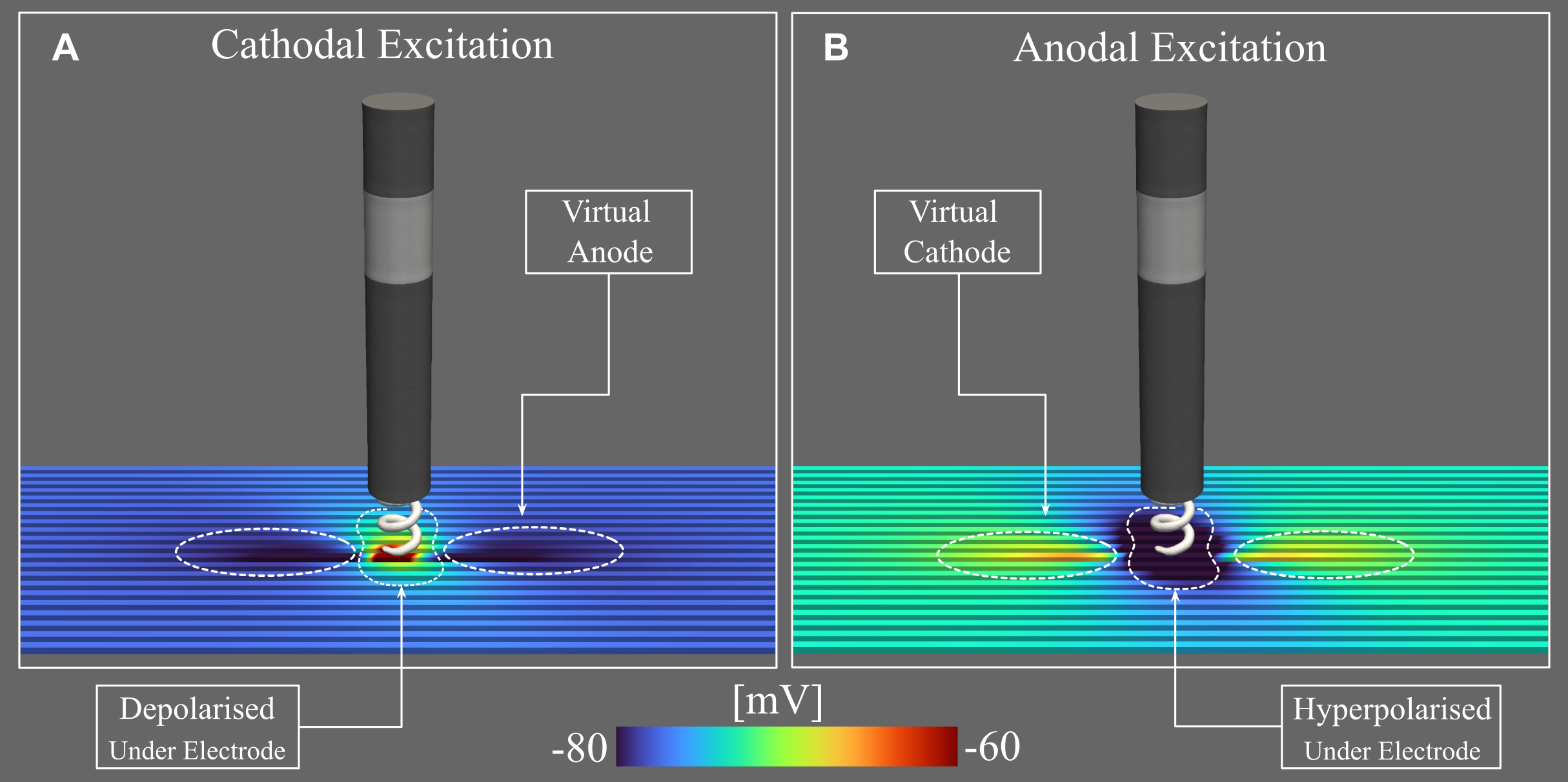}
    \caption{\Glspl{vep} of the \gls{lbb} under cathodal and anodal pacing. 
     A) Cathodal pacing resulted in direct depolarization of \gls{lbb} bundles adjacent to the tip,
     with two virtual anodes hyperpolarizing bundles longitudinally.
     B) Anodal pacing caused \glspl{vep} of opposite polarity. 
     In both cases wavefront propagation is initiated only within adjacent bundles, 
     but not across the band as a whole.}
    \label{fig:VEPs_Cathodal_vs_Anodal}
\end{figure}

\begin{figure}[H]
    \centering
    \includegraphics[width=1.0\linewidth]{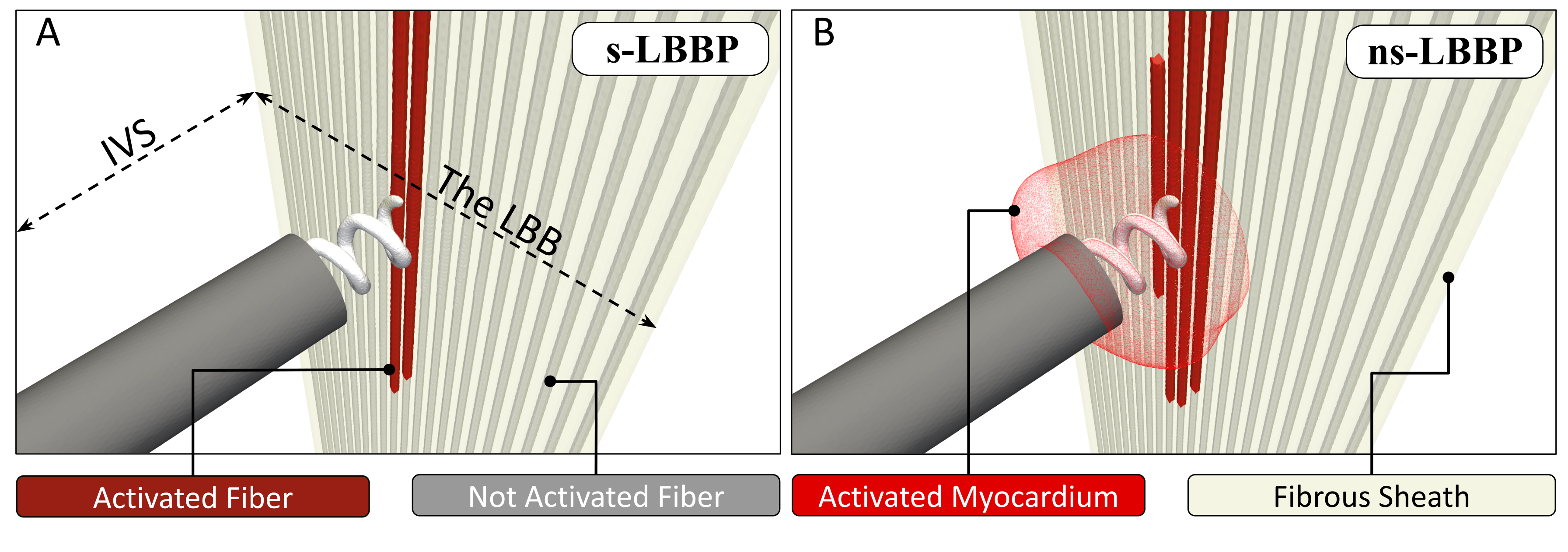}
    \caption{\Gls{lbb} capture for \gls{slbbp} and \gls{nslbbp}. 
    (A) \gls{slbbp} is achieved for tip positions touching or penetrating the \gls{lbb} band at relatively low stimulus strengths.
    (B) Increasing the stimulus strength leads to loss of \gls{slbbp} capture and transition to \gls{nslbbp}, due to additional recruitment of local myocardial tissue.
    In both cases, activation of the \gls{lbb} band remained partial, capturing only bundles adjacent to the electrode. Full capture of the entire \gls{lbb} band, enabling simultaneous activation of all fascicles, was not observed within the minimum thresholds indicated by the strength-duration curve for \gls{nslbbp}.}
    \label{fig:s-LBBP_vs_ns-LBBP}
\end{figure}

Achieving \gls{slbbp} requires the helical tip to be advanced to the \gls{lv} subendocardium 
to physically touch or penetrate the \gls{lbb} band.
\Gls{lbb} capture was highly localized for both cathodal and anodal stimulation, 
with virtual electrodes polarizing bundles largely longitudinally, 
with a very limited transverse polarization (Fig.~\ref{fig:VEPs_Cathodal_vs_Anodal}A and B). 
As such, \gls{slbbp} capture of the \gls{lbb} over the entire width of the band, 
thus ascertaining a synchronous activation of all fascicles, could not be achieved, 
only \gls{lbb} fibers in the immediate vicinity of the tip were activated (Fig.~\ref{fig:s-LBBP_vs_ns-LBBP}A).
Even for much higher strengths leading to \gls{nslbbp}, 
more bundles were simultaneously activated, but a full capture of the \gls{lbb} band was not observed within the minimum thresholds indicated by the \gls{sd} curve for \gls{nslbbp} (Fig.~\ref{fig:s-LBBP_vs_ns-LBBP}A).

% \todo[inline]{Maybe we add some more data here. At the maximum strength tested, how many fibers activated?}

\begin{figure}[H]
    \centering
    \includegraphics[width=1\linewidth]{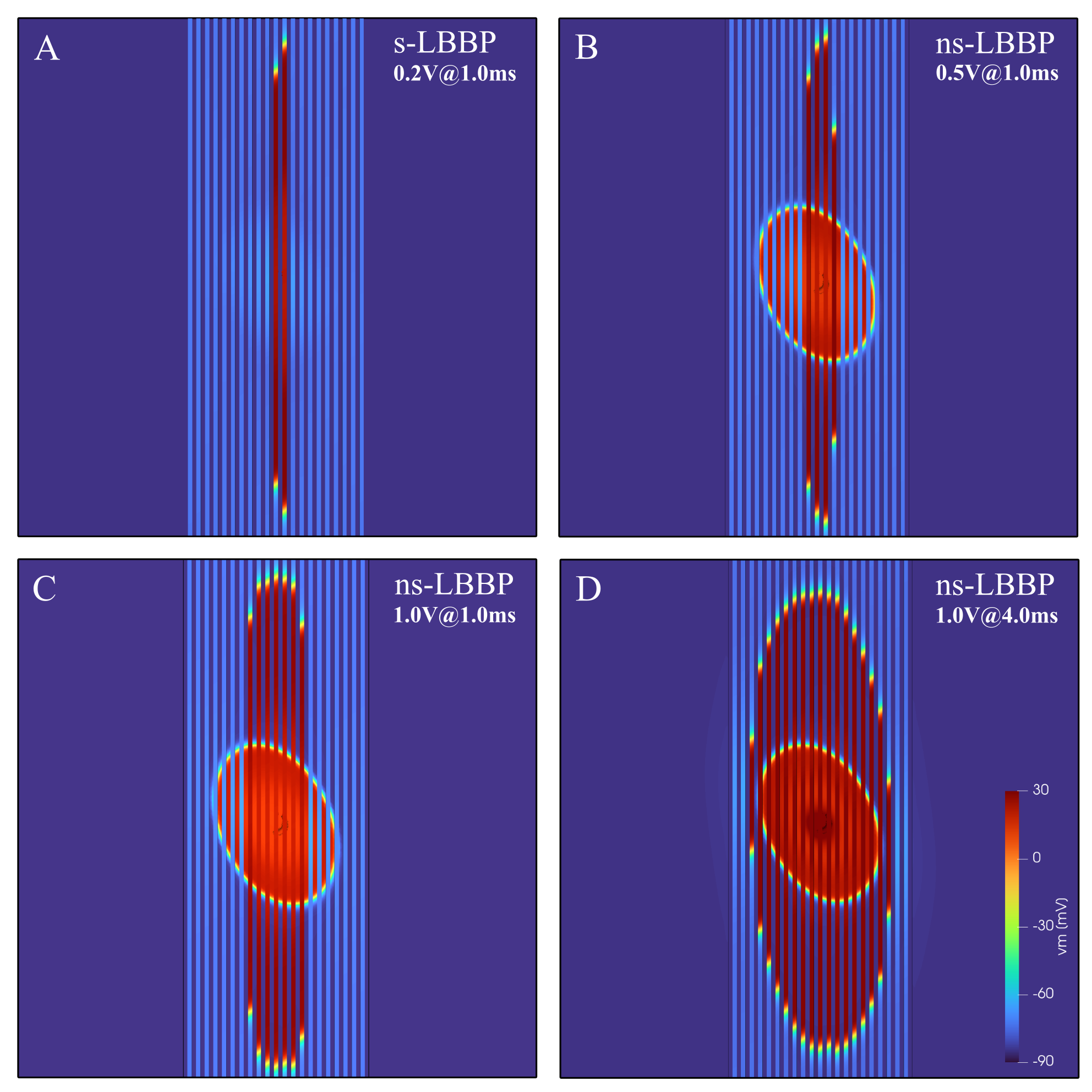}
    \caption{Differential capture of \gls{lbb} and \gls{ivs} myocardium. 
    A) Selective capture of the \gls{lbb}, with only a minor sub-threshold 
    depolarization of the \gls{ivs}. B) With increasing strength selective capture
    of the \gls{lbb} is lost. Further increasing strength C) and duration D) 
    improved \gls{lbb} capture, without affecting \gls{ivs} capture.}
    \label{fig:FibersRecruitment}
\end{figure}

\subsection{Effect of Lead Orientation on Capture Mechanisms}

To investigate the impact of a non-orthogonal deployment on \gls{sd} curves 
a lead deployment under an angle of 20$^\circ$ 
tilted longitudinally along the \gls{lbb} in the apico-basal direction, $\alpha_{\mathrm{l}}$, 
or transversely tilted in the antero-posterior direction by $\alpha_{\mathrm{t}}$ 
was considered and compared to an orthogonal deployment (Fig.~\ref{fig:advance_stories}A and B, respectively).

\begin{figure}[H]
    \centering
    \includegraphics[width=0.75\linewidth]{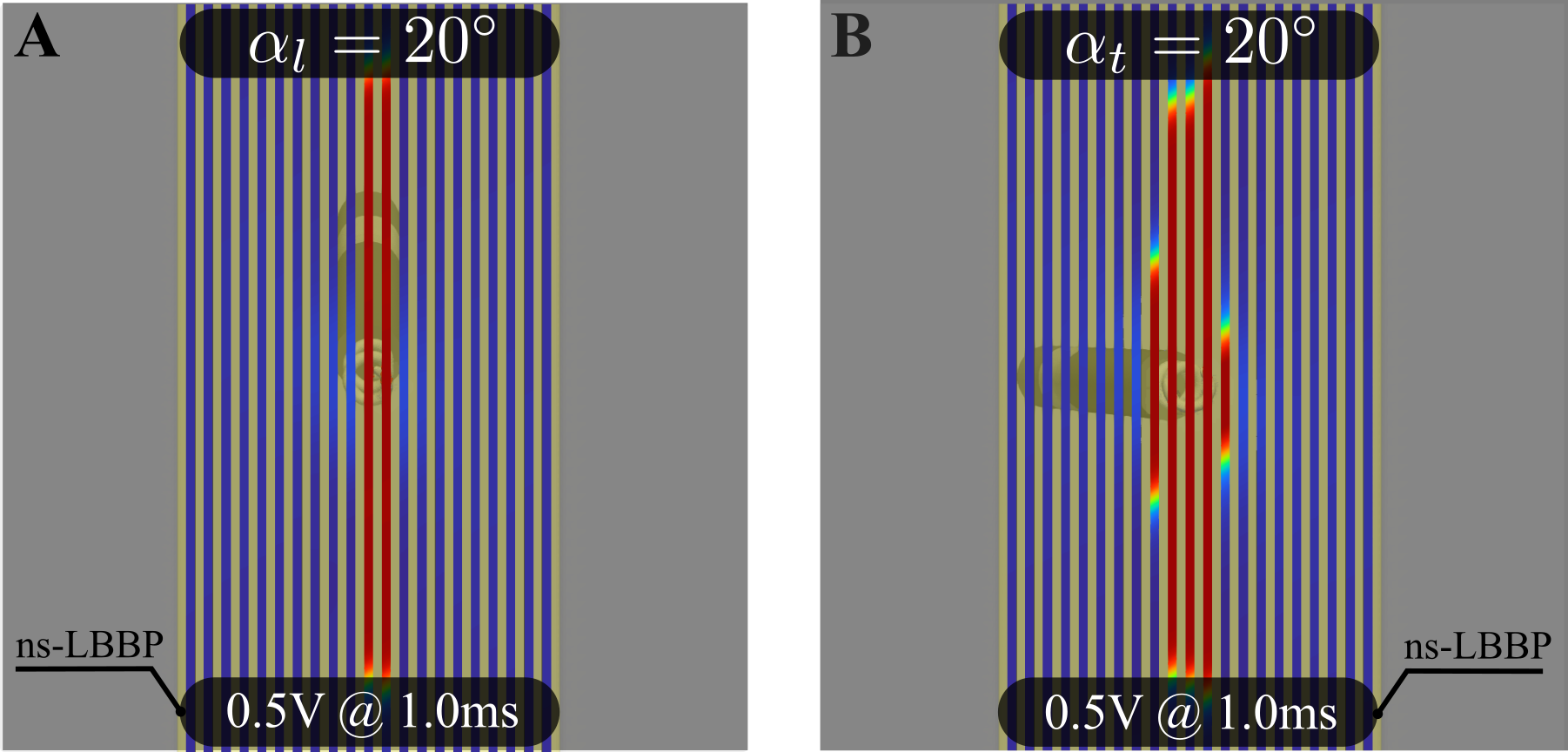}
    \caption{Effect of non-orthogonal lead orientation on \gls{lbb} fiber recruitment. Tilting the lead along (A) and transverse to (B) the \gls{lbb} axis demonstrated that transverse orientation resulted in the engagement of more number of fibers.}
    \label{fig:lead_orientation_effect}
\end{figure}

\Gls{sd} curves reveal that a longitudinal tilt by $\alpha_{\mathrm{l}}$ 
reduced the selective capture threshold at shorter pacing durations, 
from \SI{0.27}{\volt} to \SI{0.25}{\volt} at \SI{0.5}{\milli \second}, 
and effectively increased the margin to \gls{ivs} capture over all tested durations (Fig.~\ref{fig:SD_Curves_2}A). 
Conversely, a transverse tilt by $\alpha_{\mathrm{t}}$ increased the \gls{lbb} capture threshold,
from \SI{0.27}{\volt} to \SI{0.30}{\volt} at \SI{0.5}{\milli \second},
thereby narrowing the margin for \gls{slbbp} (Fig.~\ref{fig:SD_Curves_2}C).
While worse in terms of selectivity, bundle recruitment turned out to be more effective 
with a transverse tilt as the number of activated bundles increased, 
particularly for shorter duration stimuli $\le$\SI{2}{\milli \second} (Fig.~\ref{fig:lead_orientation_effect}).

\subsection{Septal Perforation}
Perforation led to a notable increase in capture threshold,
from \SI{0.27}{\volt} to nearly \SI{2}{\volt} at \SI{0.5}{\milli \second}.
For shorter stimulus duration the margin for \gls{slbbp} was significantly reduced,
becoming zero at \SI{0.5}{\milli \second} (Fig.~\ref{fig:SD_Curves_2}D). 
The opposite effect was witnessed for longer stimulus duration of $\ge$\SI{1}{\milli \second}
where the margin for \gls{slbbp} was substantially increased, up to $\approx$\SI{0.3}{\volt}.
This upward shift of the \gls{sd} curves towards higher stimulation strength 
can be attributed to the disruption of direct contact of the helical tip with the \gls{lbb}, 
along with the insulating properties of the fibrous sheath wrapped around the \gls{lbb}.
%which hinder effective current transmission. 
To corroborate the role of the fibrous sheath as a primary factor 
responsible for the increase in \gls{slbbp} capture threshold, 
a modified setup was considered
where the low conductive fibrous tissue facing the \gls{lv} endocardium was removed,
thus exposing the \gls{lbb} directly to the \gls{lv} blood pool. 
This led to a marked reduction in the \gls{slbbp} capture threshold, 
facilitating activation of the \gls{lbb} with a stimulus of only \SI{0.14}{\volt} at \SI{0.5}{\milli \second}.

\subsection{Dislodgement Effect}
\begin{figure}[H]
    \centering
    \includegraphics[width=1.0\linewidth]{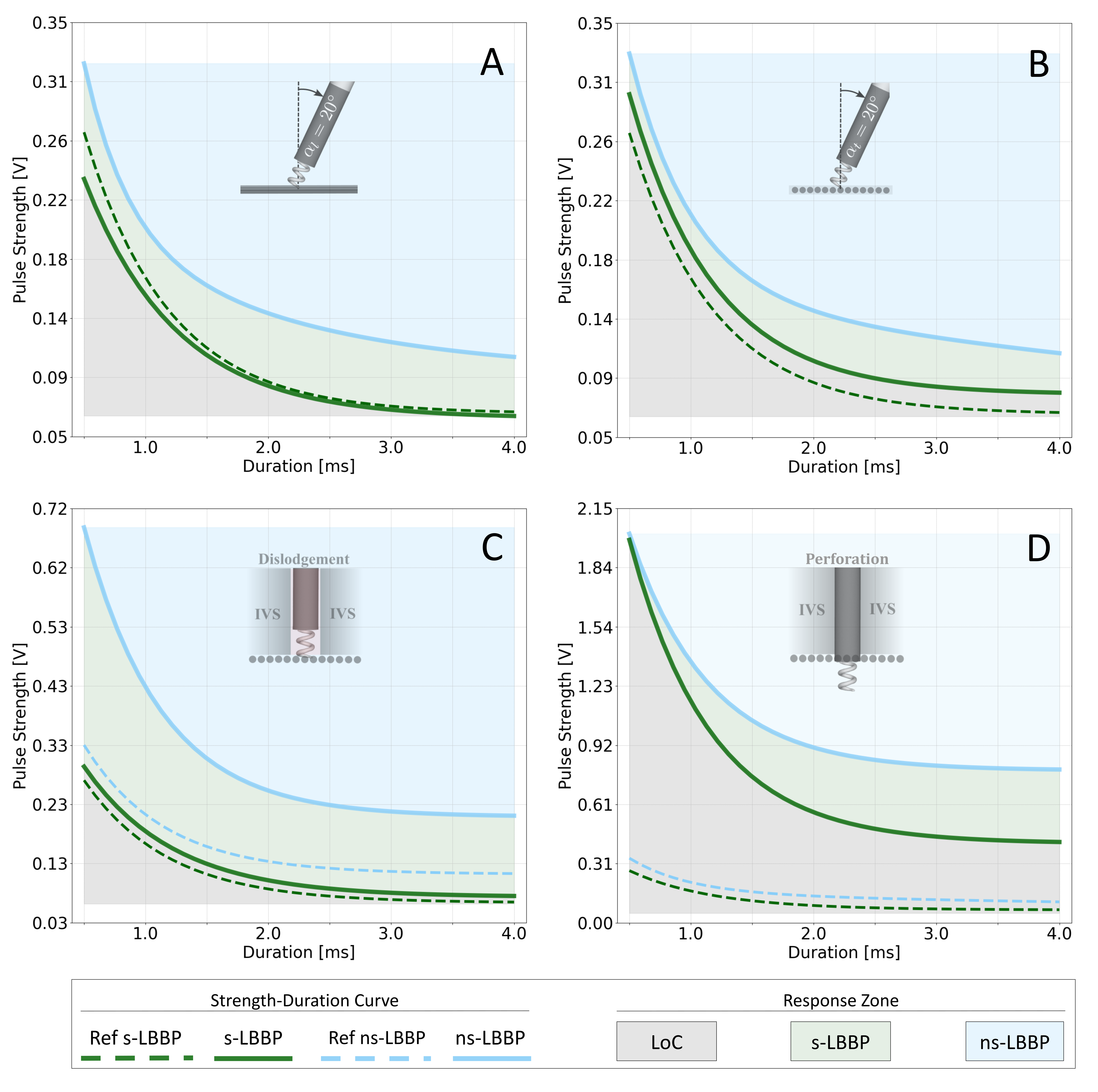}
    \caption{Alteration of \gls{lbbap} \gls{sd} curves for 
    (A) deployment of the lead with an tilt angle $\alpha_{\mathrm{l}}$
    shifted the \gls{slbbp} towards lower thresholds, 
    and expanded the margin for \gls{slbbp};
    (B) deployment of the lead with an antero-posterior transverse tilt angle $\alpha_{\mathrm{t}}$
    led to a shift of the \gls{slbbp} thresholds towards higher strengths, 
    and narrowed the margin for \gls{slbbp};
    (C) dislodgement shifted \gls{sd} curves towards higher stimulus strengths, 
    with a marked expansion of the margin for \gls{slbbp};
    (D) Septal perforation increased both \gls{slbbp} and \gls{nslbbp} capture thresholds, 
    increasing the margin for \gls{slbbp} at longer stimulus duration.
    }
    \label{fig:SD_Curves_2}
\end{figure}

\Gls{sd} curves under dislodgement conditions shifted upwards towards higher strengths,
albeit at different rates, thus notably expanding the margin for achieving \gls{slbbp}.
For instance, at short stimuli lasting \SI{0.5}{\milli \second}, 
the \gls{lbb} capture threshold increased from \SI{0.27}{\volt} observed in the reference setting to \SI{0.3}{\volt} (Fig.~\ref{fig:SD_Curves_2}C),
whereas a much more substantial increase was witnessed for capturing the \gls{ivs}, 
from \SI{0.33}{\volt} to \SI{0.7}{\volt}.

\subsection{Monitoring Impedance Trend during Implant}

\begin{figure}[H]
    \centering
    \includegraphics[width=0.5\linewidth]{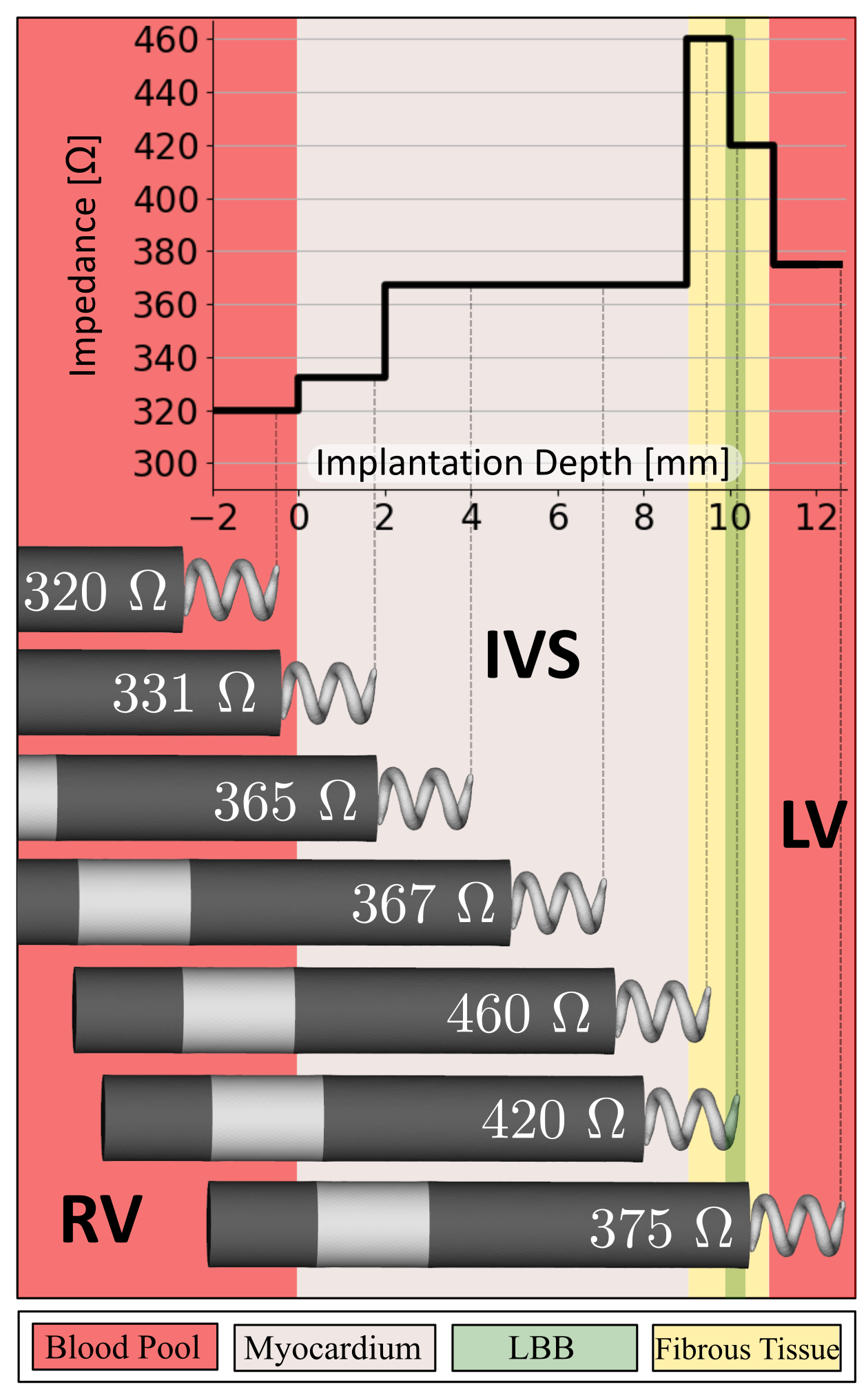}
    \caption{Variation of impedance $Z$ at different implantation depths $z$ 
    in the \gls{ivs} during \gls{lbbap}. Major changes occurred as the helix penetrated the \gls{ivs}, reached the fibrous sheath, contacted the \gls{lbb}, and perforated the \gls{ivs}.}
    \label{fig:impedance_depth}
\end{figure}

The implant procedure was simulated by computing the unipolar tip impedance $Z$
at various depths $z$, representing discrete positions during lead advancement.
Most pronounced changes in $Z$ were observed during the initial implantation phase 
where the helix gradually penetrated into the \gls{rv} subendocardium of the \gls{ivs},
and later in the \gls{lv} subendocardium 
when the tip approached the \gls{lbb} (Fig.~\ref{fig:impedance_depth}).
The observed trends as well as the computed impedances are consistent with clinical reports \cite{vijayaraman_his-bundle_2019, liu_left_2021, ponnusamy_how_2021, ponnusamy_electrophysiological_2022, ravi_late-onset_2020, ghosh_septal_2024}.
A sudden increase in $Z$ was witnessed when the tip penetrated into the low conductive fibrous sheath,
followed by a sudden drop in $Z$, indicating that the tip advanced into the \gls{lbb} and, 
thus, approached an ideal position for \gls{slbbp}. 
A final drop in $Z$ is noticed where the tip starts to perforate the \gls{lv} endocardium,
getting in physical contact with the \gls{lv} blood pool.

\subsection{Impact of Tissue Composition upon Impedance}
\begin{figure}[H]
    \centering
    \includegraphics[width=1.0\linewidth]{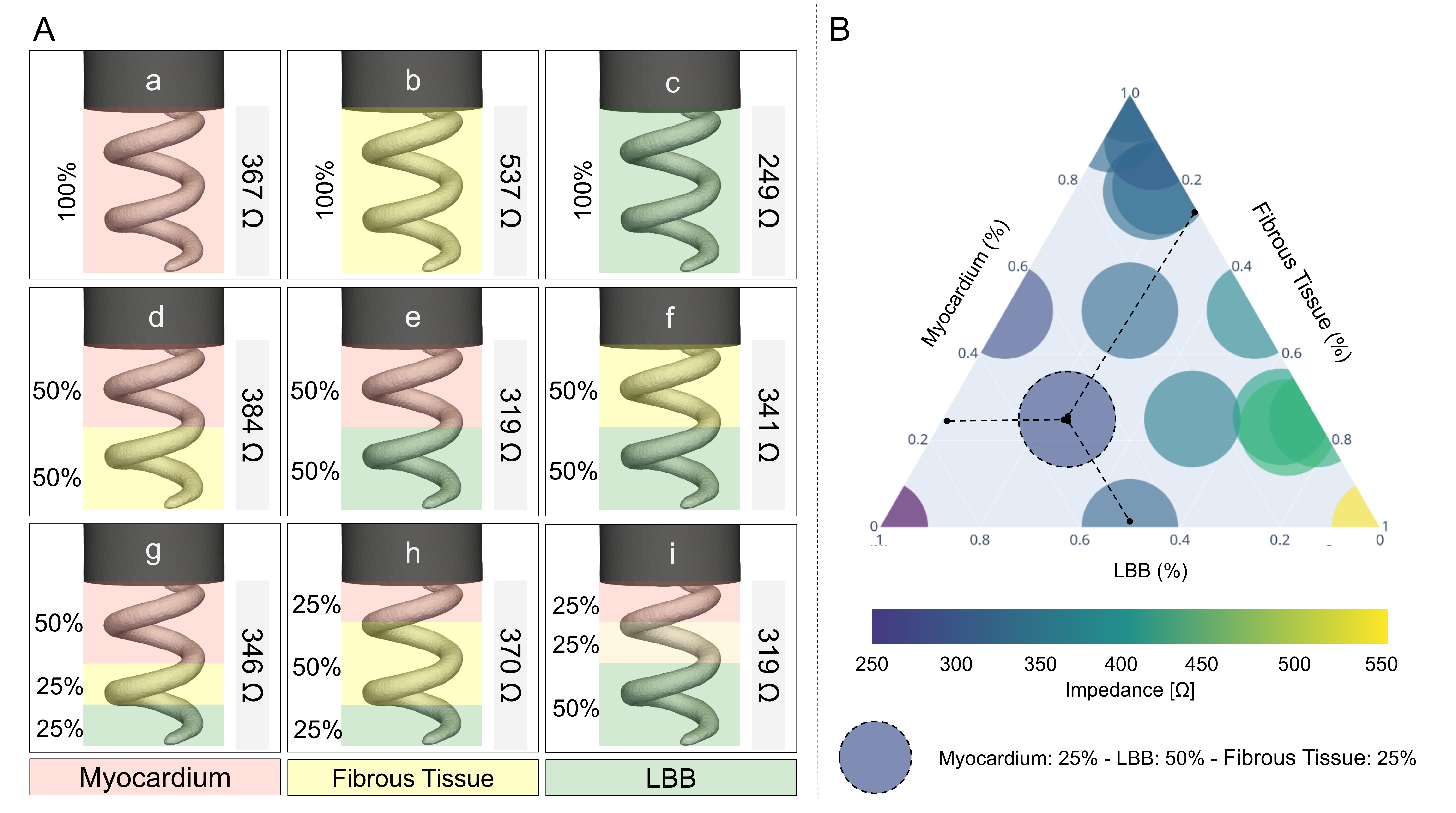}
    \caption{Analyzing the role of volumetric tissue composition at the tissue-electrode interface 
    on unipolar tip impedance. 
    A) Variation of tissue composition with estimated impedance. 
    B) Ternary plot illustrating the relationship between tissue composition and measured impedance, showing a strong positive correlation with low-conductivity fibrous tissue and a strong negative correlation with high-conductivity \gls{lbb}. Intermediate-conductivity \gls{ivs} myocardium exhibited only a weak correlation.}
    \label{fig:impedance}
\end{figure}

The unipolar tip impedance $Z$ was determined for variable volumetric tissue compositions around the tip (Fig.~\ref{fig:impedance}A), 
yielding $Z$ values ranging from 249$\Omega$ to 537$\Omega$.
For the baseline configuration used throughout this study, 
where the helix was exposed to 75\% myocardium, 25\% fibrous tissue, and 5\% \gls{lbb}, 
a $Z$ of 460 \(\Omega\) was estimated.
With the assumed conductivities our analysis indicates, as expected, 
strongly positive and negative correlations with values of $0.85$ and $-0.8$ 
for low-conductivity fibrous tissue and high-conductivity \gls{lbb}, respectively,
and a low correlation of $0.1$ for the intermediate-conductive myocardial tissue (Fig.~\ref{fig:impedance}B).
Our analysis suggests that $Z$ is largely governed by the exposure to the fibrous tissue and \gls{lbb},
whereas the \gls{ivs} myocardial tissue plays the least important role in modulating $Z$.

\subsection{Impact of Tissue Composition upon Capture Thresholds}
Analogously, the dependence of capture thresholds for \gls{slbbp} and \gls{nslbbp} 
upon volumetric tissue composition was analyzed by measuring the change in \gls{sd} curves
relative to the baseline reference configuration, 
where the helical tip was embedded in 70\% myocardium, 25\% fibrous tissue, and 5\% \gls{lbb}.

\begin{figure}[H]
    \centering
    \includegraphics[width=0.9\linewidth]{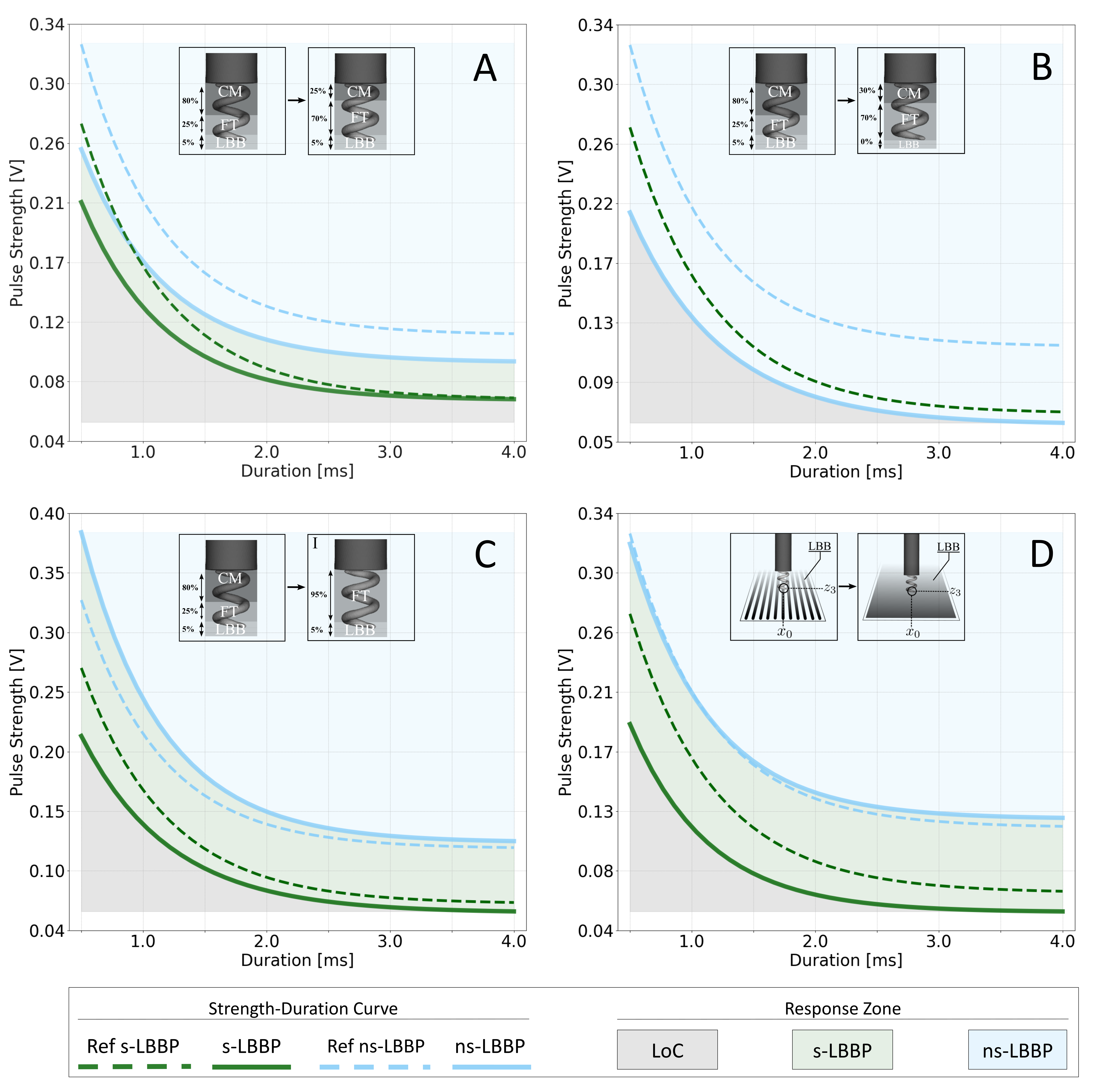}
    \caption{Dependence of \gls{slbbp} and \gls{nslbbp} thresholds upon tissue composition. 
    (A) Increasing the width of the fibrous tissue (FT) layer and maintaining contact with \gls{ivs} myocardium (CM)
    and \gls{lbb}.
    (B) Increasing FT layer towards \gls{lbb} eliminating direct physical contact with the \gls{lbb},
    and, (C) towards the \gls{ivs} eliminating physical contact of the tip with the \gls{ivs}.
    (D) Effect of tight packing of \gls{lbb} bundles upon \gls{sd} curves.}
    \label{fig:SD_Curves_3}
\end{figure}

An increase in fibrous tissue proportion to 75\%, with reduced myocardium to 20\% and constant \gls{lbb} at 5\%, resulted in a downward shift of the \gls{sd} curves for both \gls{slbbp} and \gls{nslbbp}, accompanied by a narrowing of the \gls{slbbp} margin (Fig.~\ref{fig:SD_Curves_3}A). 

Further increasing the width of the fibrous tissue layer towards the \gls{lbb} 
to eliminate the direct physical contact of the helix with the \gls{lbb} resulted in a loss of \gls{slbbp}
and a substantial reduction of the \gls{nslbbp} threshold (Fig.~\ref{fig:SD_Curves_3}B).

Conversely, widening the fibrous tissue layer towards the \gls{ivs} myocardium to isolate the tip from direct myocardial contact, while maintaining optimal helix contact with the \gls{lbb}, caused the \gls{slbbp} and \gls{nslbbp} \gls{sd} curves to shift in opposite directions, markedly enlarging the \gls{slbbp} margin (Fig.~\ref{fig:SD_Curves_3}C).

Such a configuration where the \gls{lbb} remains structurally intact and unaffected by underlying pathology, 
but is fully covered by fibrous tissue of sufficient width, can be considered optimal for achieving \gls{slbbp}. 

% Importantly, this can be exploited for lead design!

Finally, the influence of the \gls{lbb} structural configuration was investigated.
Assuming a tighter packing of bundles to form a ribbon-like structure, 
without interspersing fibrous tissue between the bundles lowered the \gls{slbbp} threshold
and increased the \gls{nslbbp} threshold, also yielding an increased margin for achieving \gls{slbbp} (Fig.~\ref{fig:SD_Curves_3}D).

\begin{figure}[H]
    \centering
    \includegraphics[width=0.75\linewidth]{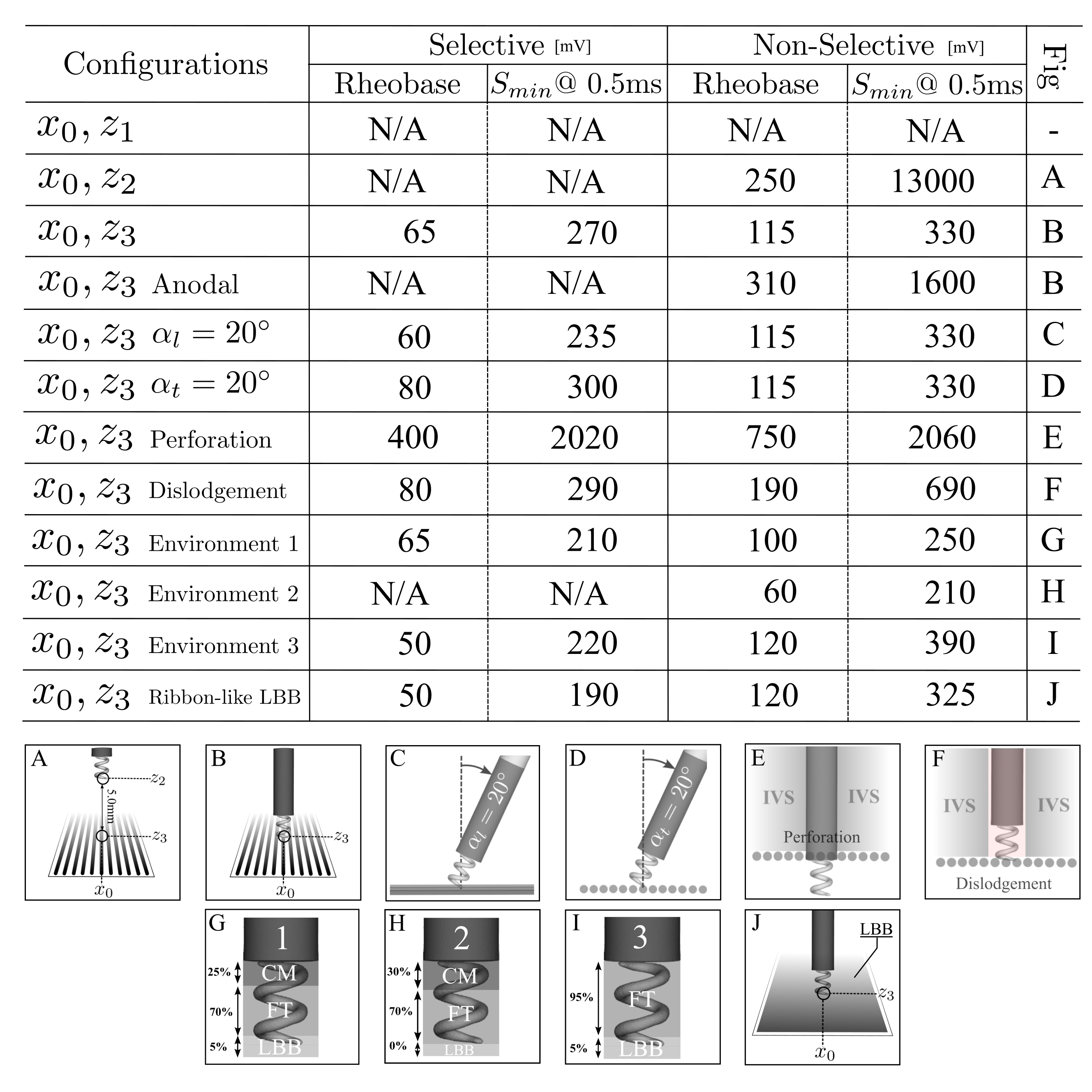}
    \caption{Rheobase and stimulus strength for \gls{slbbp} and \gls{nslbbp} capture for different conditions: 
    (A-B) variable distance of tip to \gls{lbb};
    (C-D) lead deployment tilt in long and transverse axis;
    (E) \Gls{lv} myocardial perforation;
    (F) lead dislodgement;
    (G-I) variable tissue composition around tip; and,
    (J) \gls{lbb} bundle tight packing density.
    CM: Cardiomyocyte and 
    FT: Fibrous Tissue. 
    }
    \label{fig:summary_table}
\end{figure}

% ============================================== Discussion =================================================

\section{Discussion}\label{sec:discussion}
%
% Intro discussion Ed, adjust to our study
The application scope of \gls{csp} therapy is increasingly broadening, 
as it offers potential to closely restore a normal ventricular activation sequence.
However, \gls{csp} is a complex therapy 
and fully exploiting its benefits requires optimizing lead design, deployment and pacing protocol.
Achieving optimization is hampered by an insufficient understanding of the physical mechanisms 
governing selective capture of the \gls{lbb},
and by limited insight into acute and long-term consequences of the type of capture, \gls{slbbp} versus \gls{nslbbp}.
Due to a large number of confounding factors including 
lead geometry and pacing vector, deployment location, \gls{lbb} and \gls{ivs} anatomy 
and their electrophysiological differences, and the electrical coupling conditions 
of the \gls{lbb} within the \gls{ivs}, 
developing a more comprehensive mechanistic perspective by purely experimental or clinical inquiries 
has proven elusive
since an accurate mapping of the pacing-induced activation patterns 
or differentiating \gls{lbb} from the \gls{ivs} myocardial response is challenging with current technologies.

In this study, we address these limitations by performing the anatomically and biophysically most detailed computer simulations 
of \gls{lbb} and \gls{ivs} active response to the electric field 
generated by a geometrically accurate model of a clinically widely used \gls{csp} lead,  
to investigate capture mechanisms  
to elucidate the factors determining whether activation is selective or not. 
% Now we summarize our key insights
Our study shows that 
1) robust \gls{slbbp}, 
with a sufficient margin in capture thresholds between \gls{lbb} and \gls{ivs} myocardium,
can only be achieved by placing the electrode tip in the immediate vicinity of the \gls{lbb}, 
or, by contacting or penetrating the \gls{lbb},
2) for \gls{slbbp} only a partial capture of the \gls{lbb} was witnessed,
suggesting that \gls{slbbp} may yield only \gls{lv} activation sequences similar to those produced by a fascicular block,
but that a normal healthy \gls{lv} activation cannot be fully restored,
3) a synchronous capture of all bundles within the \gls{lbb} could only be achieved for higher stimulus strengths,
but all of these were associated with \gls{nslbbp},
4) switching tip polarity to anodal pacing mode did not result in superior selectivity for capturing the \gls{lbb}, 
in contrast to previous simulation studies \cite{vigmond_how_2021}, and
5) capture thresholds of the \gls{lbb} were only moderately affected by the variability in lead deployment 
which cannot be concisely controlled during implant.

The validity of the model is supported by both quantitative and qualitative predictions 
that are consistent with clinical reports.
Predicted unipolar tip impedances \cite{vijayaraman_his-bundle_2019, liu_left_2021, ponnusamy_how_2021, ponnusamy_electrophysiological_2022} and capture thresholds \cite{burri_ehra_2023, liu_left_2021, vijayaraman_left_2019} all fell within the reported ranges,
impedance trends observed for advancing the lead towards optimal placements within the \gls{lbb} \cite{orlov_left_2023}
or for the implant complications septal perforation and lead dislodgement 
were consistent with clinical observations \cite{ghosh_septal_2024, ponnusamy_late_2021}.
This close agreement with clinical recordings suggest 
that anatomically accurate and biophysically detailed models as used in this study 
may offer mechanistic insights that can be used as feedback 
for guiding the development of optimized lead designs, pacing protocols or implant procedures.

\subsection{Capture Thresholds}
Thresholds for capturing \gls{lbb} and \gls{ivs} were investigated as a function of 
lead placement and orientation within the \gls{ivs} relative to the orientation of the \gls{lbb}, 
stimulus polarity, the tissue composition surrounding the tip as well as the impact of the implant complications 
\gls{lv} perforation and lead dislodgement. 
Parameter ranges were identified where the \gls{lbb} could be captured at a stimulus threshold 
that was well below the threshold of the \gls{ivs} myocardium, with a sufficient margin to robustly achieve \gls{slbbp}.
These regions were limited to deep implant locations at $z_3$ where the tip was in physical contact with the \gls{lbb}.
Differences in capture thresholds between \gls{ivs} myocardium and \gls{lbb} were rather small,
with the largest observed difference of \SI{95}{\milli \volt}.
Importantly, with the given device design only a partial capture of the \gls{lbb} was witnessed 
for the stimulus strength leading to \gls{slbbp},
inducing \gls{lv} activation sequences resembling a fascicular block, but not a normal healthy activation.

\subsubsection{Lead Placement Depth within the \glsentryshort{ivs}}
Key to achieving \gls{slbbp} was a deep implant of the lead into the \gls{ivs}
such that the helical tip was in very close proximity or in contact with the \gls{lbb}. 
With increasing distance between electrode tip and \gls{lbb} the threshold for capturing the \gls{lbb} rose quickly.
For a midwall tip placement at $z_2$, \gls{lbb} capture could be achieved, 
but the threshold for \gls{nslbbp} rose markedly, 
from \SI{0.33}{\volt} at the deep implant location $z_3$ to the \SI{12}{\volt} at $z_2$.
Thus, lead placement in the \gls{ivs} at a distance from the \gls{lbb} most likely leads to \gls{lvsp} only and not \gls{nslbbp}, 
as only the \gls{ivs} myocardium will be captured, without exciting the \gls{lbb}~\cite{moustafa_conduction_2023}.
Although \gls{lvsp} was not previously classified as \gls{csp} \cite{burri_ehra_2023}, it is now recognised as part of the \gls{csp} spectrum \cite{glikson_european_2025}. While it does not offer the same electrophysiological and hemodynamic benefits as \gls{lbbp}, it has been shown to restore \gls{lv} synchrony and improve hemodynamic function to a degree comparable with conventional \gls{crt} using \gls{bivp}~\cite{heckman_physiology_2021, curila_lvsp_2024}. 
In our model a minimum capture threshold with a rheobase of \SI{0.065}{\volt} 
or with a  stimulus of \SI{0.27}{\volt} at \SI{0.5}{\milli \second} was identified for \gls{lvsp},
which is in line with reported values~\cite{wang_efficacy_2020, burri_ehra_2023}.

The \gls{ivs} myocardial capture threshold also varied with depth, but to a marginal extent.
For a stimulus duration of \SI{0.5}{\milli \second} at implant depths $z_1$, $z_2$ and $z_3$ 
myocardial capture occurred at \SI{0.24}{\volt}, \SI{0.27}{\volt} and \SI{0.33}{\volt}, respectively,
consistent with a previous report~\cite{vigmond_how_2021}.
This progression was attributable to factors such as conductive heterogeneity between tip and ring, 
due to the presence of blood pools or fibrous tissue, 
and a reduced electrotonic load in the vicinity of the \gls{ivs} surface.

\subsubsection{Lead Deployment Angle}
Guidelines for \gls{lbbap} recommend perpendicular lead implantation into the \gls{ivs} 
for safer and easier lead fixation~\cite{de_pooter_guide_2022, shroff_comparison_2024}. 
However, in practice, the deployment angle within the \gls{ivs} cannot be precisely controlled and is not known during implant. 
Thus, the impact of deviating from a perpendicular lead implant has not been a primary focus in \gls{lbbap} studies~\cite{burri_ehra_2023, cano_left_2021, liu_left_2021, vijayaraman_left_2019}. 
In contrast, computational simulations enable a controlled variation of the deployment angle 
and the mechanistic exploration of its impact.
Mechanistically, an altered deployment angle affects the orientation of the electric field 
relative to the \gls{lbb}, 
and, thus, the activating function driving the change in tissue polarization (see in Eq.~\eqref{eq:_bidm_parab}).
For a first evaluation of its potential role an angular deviation relative to an orthogonal reference implant of $20^\circ$ 
was considered.

Tilting the lead increased the potential gradient $\nabla \phi_{\rm e}^p$ 
within a plane parallel to the \gls{lbb}.
For a longitudinal tilt this led to a lower capture threshold and enhanced selectivity 
by increasing the margin between \gls{lbb} and \gls{ivs} capture (Fig.~\ref{fig:SD_Curves_2}A).
The higher potential gradient acting along the bundles of the \gls{lbbb}
was able to create a sufficiently large cathodal polarization at a lower stimulus strength.
Tilting the lead transversely to the \gls{lbb} led to a higher threshold and reduced selectivity 
(Fig.~\ref{fig:SD_Curves_2}B)
as the potential gradient was oriented transversely to the bundles 
which were assumed to be electrically isolated from each other. 
Thus, a higher field gradient was needed here to create a 
critical supra-threshold cathodal polarization along the bundles. 
However, the increased transverse electric field led to a more effective recruitment of fibers of the \gls{lbb}, 
suggesting that such an oblique implant may facilitate the selective capturing of the \gls{lbb} in its entirety
rather than achieving a partial capture of a few fibers only.
Clinically, this attribute may be significant for \gls{lbbb} patients
as activating additional fibers could improve \gls{lbbb} correction,
even if it comes at the cost of slightly higher \gls{nslbbp} thresholds.

\subsection{Anodal Excitation in \glsentryshort{lbbap}}

In \gls{lbbap}, anodal excitation occurs when the lead is implanted sufficiently deep into the \gls{ivs}
such that the ring electrode lead is in physical contact with the \gls{ivs}. 
Chances for anodal capture of the \gls{rv} septal wall are enhanced in situations where an increased pacing output is needed 
at an average threshold of 3.3$\pm$\SI{1.6}{\volt} at a duration of \SI{0.5}{\milli \second} \cite{hopman_loss_2023, arnold_his-purkinje_2020, tamborero_anodal_2006}.
However, to the best of our knowledge, the use of the helical tip as the anode in \gls{lbbap} has not been reported in the literature.
Reversing the polarity of the helix to function as the anode resulted in a significant increase in capture thresholds 
for both \gls{ivs} myocardium and \gls{lbb}~(Fig.~\ref{fig:SD_Curves_1}D).
Due to the higher stimulus strengths required, \gls{slbbp} could not be achieved 
over the range of tested stimulus durations up to \SI{3}{\milli \second},
only \gls{nslbbp} was witnessed.
Moreover, anodal capture of the \gls{lbb} was more susceptible to lead placement,
where placing the lead slightly off the center of the \gls{lbb} 
resulted in a marked increase in capture threshold for \gls{nslbbp} 
-- from \SI{1.58}{\volt} to \SI{8.2}{\volt} at \SI{0.5}{\milli \second} --
due to the virtual cathodes not fully encompassing the \gls{lbb} domain.

The markedly higher anodal capture thresholds are explained by fundamental differences in the activation mechanism. 
Unlike cathodal excitation, where depolarisation occurred directly in tissue adjacent to the physical electrode, 
anodal stimulation initially induced hyperpolarisation there \cite{wikswo_virtual_1995}~(Fig.~\ref{fig:VEPs_Cathodal_vs_Anodal}A-B). 
This hyperpolarised region was surrounded by virtual cathodes, which facilitated depolarisation in nearby tissue regions (Fig.~\ref{fig:VEPs_Cathodal_vs_Anodal}).
The excitation mechanism was either anodal-make or break, depending on the stimulus strength~\cite{sambelashvili_virtual_2004}.
For stronger stimuli, the virtual cathodes may polarize sufficiently to achieve anode-make excitation, 
for the weaker stimuli more relevant for achieving \gls{slbbp} anode-break excitation was the dominant mechanism.

Our results do not corroborate the finding of better selectivity as reported previously~\cite{vigmond_how_2021}.
This discrepancy may be explained by differences in lead design and experimental parameters. 
In~\cite{vigmond_how_2021}, a wider range of conceivable lead implant depths and orientations was evaluated,
including deep lead implant oriented parallel to the \gls{lbb}, 
only a single fiber without electrically insulating fibrous tissue was considered,
and enhanced selectivity was observed for longer stimuli lasting $>$\SI{4}{\milli \second}.
In our study deployment angles as well as stimulus duration were limited 
by the bounds of technical feasibility imposed by the given lead design. 
Our physically most detailed simulations do not support the notion of enhanced selectivity with anodal stimulation.
Rather, cathodal stimulation appears more effective, due to the significantly lower capture threshold.

\subsection{Implant Complications in \glsentryshort{lbbap}}

Dislodgement may occur when drilling fails to adequately penetrate the \gls{ivs}, leading to the formation of a fluid-filled tunnel caused by local oedema~\cite{liu_left_2021}. Given its potential long-term consequences, including a reduction in \gls{lv} ejection fraction~\cite{ponnusamy_late_2021}, careful monitoring during both implantation and follow-up is recommended to minimise the risk of adverse outcomes~\cite{vijayaraman_left_2022}.

Our study indicates a modest increase in  capture threshold for \gls{slbbp}, 
from \SI{0.27}{\volt} to \SI{0.29}{\volt} at a pulse duration of \SI{0.5}{\milli \second}~(Fig.~\ref{fig:SD_Curves_2}C). 
While, quantitatively, the increase in capture threshold is thought to be dependent on the extent of dislodgement~\cite{burri_ehra_2023}, 
categorized as a macro- or a micro-dislodgement as modeled in this study,
our results are in agreement with a clinical observations, 
reporting an increase in the \gls{slbbp} threshold ~\cite{ponnusamy_late_2021}.

Furthermore, \gls{nslbbp} thresholds also rose substantially, 
from \SI{0.33}{\volt} to \SI{0.69}{\volt} at a pulse duration of \SI{0.5}{\milli \second},
potentially reflecting a loss of direct contact between electrode and \gls{ivs}.  

%\subsubsection{Septal Perforation}
Implanting the lead deep in the \gls{lv} subendocardium of the \gls{ivs} to achieve contact with the \gls{lbb}
poses the risk of septal perforation~\cite{huang_novel_2017}. 
Its identification is critical to prevent long-term complications such as 
thromboembolism and ventricular fibrillation,
requiring immediate correction by repositioning the lead~\cite{burri_ehra_2023, ponnusamy_electrophysiological_2022, vinther_late_2024}.
During perforation the helical tip is increasingly exposed to the \gls{lv} blood pool,
resulting in a substantial rise in the \gls{slbbp} capture threshold.
In our simulations a rise 
from \SI{0.27}{\volt} to \SI{2}{\volt} at a pulse duration of \SI{0.5}{\milli \second} was witnessed~(Fig.~\ref{fig:SD_Curves_2}D)
which is consistent with clinical studies 
reporting a mean capture threshold of \SI{3.02}{\volt} $\pm$ \SI{0.7}{\volt} at the same pulse duration \cite{ponnusamy_electrophysiological_2022, ravi_late-onset_2020}. 
This threshold elevation was attributable to two factors,
the loss of contact between the helix and the \gls{lbb}, leading to exposure to the \gls{lv} blood pool, 
and the presence of fibrous tissue enclosing the \gls{lbb}. 
Therefore, a sudden increase in capture threshold when the tip of the lead is approaching the \gls{lv} subendocardium,
could serve as an indicator of septal perforation.
% We removed this specific experiment

% Finally, the role of the latter was assessed by removing the fibrous tissue insulating the \gls{lbb} 
% from the highly conductive \gls{lv} blood pool.
% This led to the opposite effect, a marked decrease in capture threshold was observed,
% enabling \gls{ivs} and \gls{lbb} activation at a low threshold of \SI{0.17}{\volt} at \SI{0.5}{\milli \second}. 

% Beyond the scope, removed
% \subsubsection{Lead Fracture}

% off topic here, commented out:
% Lead fracture, while rare, represents a potential complication of CSP. It is characterized by a sudden increase in pacing impedance (typically exceeding 3000 $\Omega$), often accompanied by a complete loss of capture. Although not directly investigated in this study, it remains an important consideration in the clinical management of pacing therapies \cite{castagno_is_2024, rangaswamy_late_2024}.

\subsection{Impedance Trend Monitoring}

The absolute values measured for $Z$ strongly depend on the volumetric composition of the tissue
surrounding the lead between tip and ring electrode.
For implant depths around the optimal placement for achieving \gls{slbbp} the composition 
can be assumed to be highly variable between individuals 
due to variability in \gls{lbb} anatomy and the distribution of fibrous tissue enclosing it.
This is in line with clinical studies reporting a wide range in measured impedances \cite{wang_lbbap_2021}.
However, while absolute measures of $Z$ will inevitably vary, 
the relative trends observed during lead advancement should be robust 
and largely independent of the exact tissue composition or their conductivities. 
Advancing the tip into low-conductivity fibrous tissue will increase $Z$
whereas $Z$ will drop if the tip is further advanced into the \gls{lbb}.
Thus, in \gls{lbbap} impedance monitoring combined with fluoroscopy was suggested 
as general guidance for implanting the lead at the appropriate location, 
and was also recommended to be documented~\cite{vijayaraman_prospective_2019}. 
A rise in $Z$ serves as a marker indicating that the pacing tip was positioned close enough to the \gls{lbb} 
to facilitate \gls{slbbp} \cite{vijayaraman_prospective_2019, orlov_left_2023}. 

We aimed to quantify this trend in $Z$ by simulating lead advancement,
and measured $Z$ at selected implant depths $z$.
At locations where the tip was in contact with \gls{ivs} myocardium only, 
$Z$ values of 360 $\pm8$ $\Omega$ were recorded.
For deep implant depths where the tip was in physical contact with \gls{lbb} and fibrous tissue,
regardless of whether the response type was \gls{slbbp} or \gls{nslbbp}, 
$Z$ increased up to 460 $\pm20$ $\Omega$.
The predicted impedances were in the lower range of clinical reports,
ranging from 400 $\Omega$ up to 800 $\Omega$ \cite{huang_beginners_2019, vijayaraman_his-bundle_2019, wang_efficacy_2020, ravi_late-onset_2020}.
This discrepancy may be attributed to modeling assumptions such as representing the \gls{ivs} 
as anatomically and structurally normal, with healthy electrophysiology, 
whereas for patients requiring \gls{lbbap} this may not be the case~\cite{vijayaraman_deep_2020, vijayaraman_prospective_2019, diaz_emerging_2023}.

The simulated impedance values were governed by the variable conductivities of the tissue
in which electrode tip and ring were embedded in.
These are the highly conductive \gls{lbb},
the low-conductivity fibrous tissue, and the myocardial tissue of the \gls{ivs}
exhibiting an intermediate conductivity (see~\ref{sec:method}.2). 
Systematically varying the coverage of the helix by these tissue types revealed 
their specific influence upon impedance values. 
Covering the tip with \gls{lbb} or fibrous tissue resulted in a reduction or increase in $Z$, respectively. 
Results of this impedance analysis suggest intra-procedural measurements of $Z$ may be used 
for estimating the tissue type surrounding the helix.
Using more specific lead designs the accuracy of inferring the tissue from $Z$ may be significantly improved, 
thus providing the operator with accurate information of the tip position within the \gls{ivs}.

% Most of the following we have discussed above under complications
Finally, changes in $Z$ during lead advancement serve as an indicator of complications.
%The impedance analysis can also warn of potential complications associated with LBBAP, such as septal perforation and lead dislodgement. 

In the case of septal perforation, our model predicted a sudden drop in impedance from 460 $\Omega$ to 375 $\Omega$, which is consistent with reported trends in the literature. However, published studies typically report slightly larger reductions exceeding 200 $\Omega$. This discrepancy may be attributed to the assumption of a structurally and electrophysiologically healthy \gls{ivs} in our model.

A similar drop in $Z$ was witnessed for lead dislodgement, from 460 $\Omega$ to 380 $\Omega$, 
exhibiting behavior akin to that observed clinically. 
As such, impedance measurements alone cannot be used to discriminate these different conditions,
but combined with capture threshold data and real-time fluoroscopic images a discrimination appears feasible. 
Notably, in septal perforation, the \gls{slbbp} capture threshold rises more significantly than for dislodgement, thus potentially providing a robust diagnostic marker.

% \todo[inline]{continue here}

\subsection{Model Assumptions related to \glsentryshort{ivs} and \glsentryshort{lbb}}

The anatomical and electrophysiological modeling of the \gls{ivs} and \gls{lbb} in this study builds on numerous assumptions,
that are informed by literature reports.
The anatomy of the proximal \gls{lbb} is typically described as a ribbon-like structure
that originates just below the membranous septum, positioned between the non-coronary and right coronary aortic cusps. 
It extends distally and slightly anteriorly for approximately 10 to \SI{15}{\milli \meter} 
before branching into its fascicles \cite{titus_normal_1973, uhley_visualization_1959, demoulin_histopathological_1972}. 
The \gls{lbb} is composed of strands of Purkinje cells running nearly parallel distally, 
and is encased in a fibrous sheath \cite{james_fine_1971, elizari_normal_2017}. 
These strands, aligned in a single plane, form the characteristic ribbon-like structure 
behind the \gls{lv} sub-endocardium \cite{elizari_normal_2017}. 
Histopathological and morphometric studies in pig hearts report a mean thickness of $\approx$\SI{170}{\micro \meter} 
for the thick Purkinje fibers within the septal endocardial layer 
where the proximal \gls{lbb} is located \cite{garciabustos_quantitative_2017}. 

It was further assumed that the Purkinje strands at the proximal \gls{lbb} are electrically insulated from one another and therefore, transverse electrical coupling between Purkinje strands is considered to be restricted to the penetrating His bundle.

Our study is based on an anatomically simplified representation of the \gls{lbb} as an array of $20$ cables arranged in parallel to form a ribbon-like structure with a total width of \SI{\approx 10}{\milli\meter}. 

Altering the \gls{lbb} structure, as shown in Fig.~\ref{fig:SD_Curves_3}D, by assuming tightly packed fibers with transverse coupling in a different configuration resulted in a markedly lower threshold for \gls{lbb} activation, thereby increasing the margin for achieving \gls{slbbp}. This highlights the importance of modelling the \gls{lbb} with sufficient physiological accuracy to reliably capture its electrical response.

% It is known that the geometry of the proximal His bundle is highly variable between individuals [25] 
% and has been classified into three broad types [26]. Fibers come off of the atrioventricular node from various locations, 
% at various angles, and run together to form a structure which is more ribbon-like than a thick strand [16].

A further critical factor to consider is the transmural depth of the basal \gls{ivs}
which is know to be variable between individuals or due to type and progression of a disease. 
Based on echocardiographic measurements in a clinical study assessing the feasibility of \gls{lbbap} 
a range of $7-$ \SI{18}{\milli \meter} has been reported~\cite{vijayaraman_prospective_2019}.
A depth of \SI{10}{\milli \meter} trending rather towards an \gls{ivs} depth representative of healthy individuals 
facilitates an embedding of the entire lead for a deep implant. 
That is, both ring and tip of the lead are in contact with the \gls{ivs} tissue.

% Assumptions on the 
Differences in electrophysiological properties between \gls{lbb} and \gls{ivs} tissue under heart failure conditions, such as fibrotic remodelling of the \gls{ivs} or conduction slowing within the \gls{lbb}, have not been addressed.

% These structural features, derived from existing morphological datasets, 
% combined with the Stewart ionic model characterising the ionic dynamics of Purkinje cells, 
% established a physiologically relevant framework for the computational simulations performed in this study\cite{stewart_mathematical_2009}.

%All implantation depths were reported relative to the specific IVS geometry used in the model, 
%and the results related to penetration depth were inherently tied to this mesh geometry.  

%\subsection{Anatomical Modeling of \gls{lbb} and \gls{ivs}}
% Uncertainty of model predictions due to model assumptions on \gls{ivs}, \gls{lbb} and fibrose sheath geometry 
% and their respective conductivities show an impact upon impedance values, but without affecting the relative trends 
% during lead advancement that are important for optimizing lead placement during implant.

\subsection{Implications for Lead Deplyoment and Design}

The findings of this study underscore the critical role of lead positioning relative to the \gls{lbb} 
in determining the capture mechanism and response type,
and offer guidance for optimizing the design of electrodes to improve selective capturing capacity.
Our results indicate that successful \gls{slbbp} requires implanting the electrode tip 
directly adjacent to and in contact with the \gls{lbb}.
This is widely supported in the literature by the presence of the \gls{lbb} current of injury~\cite{saxonhouse_current_2005, burri_ehra_2023, shali_current_2022},
indicative of the lead tip penetrating the fibrous insulation, being in direct contact with the \gls{lbb} \cite{su_electrophysiological_2020}. 
Further reducing the exposure of the tip to myocardial \gls{ivs} tissue contributed to lowering the \gls{slbbp} capture threshold 
while increasing the \gls{nslbbp} threshold, ultimately leading to an expansion of the \gls{slbbp} capture zone.

% Subsequent impedance analyses revealed that tissue conductivity significantly influenced the recorded impedance values. 
% This insight was further clarified by varying the percentage of tissue type covering the helix 
% in the same optimal position close to the \gls{lbb}. 
% It was determined that differences in tissue conductivity impacted the capture threshold and altered the strength-duration curve.  

Thus, the margin facilitating \gls{slbbp} can be widened 
by increasing the proportion of low-conductivity tissue surrounding the helix 
while maintaining optimal contact with the \gls{lbb}.
This led to a decrease in capture thresholds 
for both \gls{slbbp} (\SI{0.27}{\volt} to \SI{0.21}{\volt}) and \gls{nslbbp} (\SI{0.33}{\volt} to \SI{0.26}{\volt}) at \SI{0.5}{\milli \second}, 
while narrowing the \gls{slbbp} capture zone. 
When the tip lost contact with the \gls{lbb} and was fully covered by fibrous tissue, 
\gls{slbbp} was eliminated, and the \gls{nslbbp} threshold dropped further to \SI{0.22}{\volt}. 
Conversely, when $5\%$ of the helix maintained LBB contact, with the remainder exposed to fibrous tissue, 
the \gls{slbbp} threshold again decreased to \SI{0.21}{\volt}, 
while the \gls{nslbbp} threshold increased to \SI{0.38}{\volt}, thus expanding the margin for \gls{slbbp}. 

%Final Summary on Impedance should be here
%\subsubsection{Further Considerations}

Another noteworthy limitation of \gls{slbbp} is its inability to achieve full capture of the \gls{lbb}, primarily due to the low stimulus strengths required for selective activation. Attempts to increase stimulus strength to recruit all \gls{lbb} fibers invariably resulted in \gls{nslbbp}, accompanied by premature activation of the \gls{ivs}. Consequently, \gls{slbbp} induces a temporal delay in the activation of some distant fibres relative to the pacing site, which may translate into slower overall \gls{lv} activation.

Furthermore, during \gls{nslbbp}, activation of the \gls{ivs} originates at a subendocardial site within the \gls{lv}, subsequently propagating toward the epicardial surface and progressing toward the \gls{rv} side. This activation sequence may therefore alter the normal activation pattern within the \gls{rv}.
Moreover, under \gls{lbbp}, the right bundle branch activates retrogradely via the \gls{lbb}, resulting in delayed \gls{rv} activation.

These observations support the notion that \gls{hbp} may offer superior ventricular resynchronization compared to \gls{lbbap}, as right ventricular activation delays can be avoided~\cite{ali_comparison_2023} 

\subsection{Limitations}

This study uses a structurally and biophysically highly detailed model of the \gls{ivs} 
and the embedding of the \gls{lbb} in it to investigate capture mechanisms. 
While all anatomical, structural and electrophysiological factors were informed by a wealth of data reported in the literature
it is important to note that the uncertainty on virtually all these data is significant.
Moreover, as our current understanding of the cardiac conduction system and its embedding within the myocardium remains incomplete,
and considering the paucity or even absence of data on various specific aspects of structure and function,
these knowledge gaps were filled based on assumptions which increase the uncertainty of model predictions. 
Key model assumptions comprise the diameter of bundles, the distance between these, 
and the absence of any transverse electrical coupling between the intracellular spaces 
which, if present, would yield substantially lower capture thresholds for the \gls{lbb}, as shown in Fig.~\ref{fig:SD_Curves_3}D. 
Finally, beyond geometric factors all model parameters related to electrical conductivities are notoriously challenging to measure, 
showing a wide variability \cite{roth2002electrical},
and the representation of active electrophysiological behavior in response to stimulation with electric fields 
may vary markedly among individuals, leading to discrepancies between model and physical reality.
While the impact of some of these has been quantified,
owing to the vast dimension of the parameter space we refrained from performing a thorough uncertainty quantification
of all model parameters.
Rather, we compare against clinical data and, where available, their distribution,
and emphasize in our assessment the comparison of relative trends rather than absolute measures.
These are more robust than absolute values that are known to be highly variable from clinical measurements~\cite{vijayaraman_prospective_2019, wang_lbbap_2021}.

% ============================================== Conclusion =================================================

\section{Conclusion}\label{sec:conclusion}

%LBBAP, within the context of CSP, has shown promising results in CRT. Given the increasing attention it has received over the past decade, coupled with its growing clinical application, this study was conducted to investigate LBBAP from a different perspective. 

In our study we used a structurally and biophysically detailed computer model of \gls{ivs} and \gls{lbb}
to quantitatively elucidate the role of lead position, orientation and polarity in achieving optimal \gls{slbbp} thresholds
for a clinically widely used \gls{csp} lead. 
A deep implant within the \gls{lv} sub-endocardium to ensure direct contact between electrode and \gls{lbb} 
is key for effective \gls{slbbp}.
Switching the tip polarity to anodal was not beneficial, requiring higher strengths to activate the \gls{lbb}. 
Lead orientation relative to the \gls{lbb} bundles was found to influence the \gls{slbbp} capture threshold 
and the number of synchronously activating bundles.
Our simulations replicate the impedance trends observed during implant, explaining the high variability
when advancing the tip through fibrous tissue and \gls{lbb},
as well as the significant drops in impedance were found associated with septal perforation and lead dislodgement
where the latter could be distinguished by a \gls{slbbp} capture threshold that remained nearly unchanged.
Consistence with clinically observed trends support model credibility and suggest 
that modeling may offer an effective approach for improving the design of tailored \gls{csp} leads,
facilitating a selective and synchronous activation of the entire \gls{lbb}.

% Overall, the results of this study aimed to reveal the CSP mechanism within the framework of LBBAP, 
% with the goal of advancing our understanding of this emerging CRT strategy to enhance patient outcomes and achieve better results. 

\pagebreak

\subsection*{Acknowledgements}
This research was supported by the Austrian Science Fund FWF grants grants no. 
10.55776/I6540 and 10.55776/I6474 to G.P.  
%To facilitate open access, the author has applied a CC BY public copyright license to any author-accepted manuscript resulting from this submission.

\subsection*{Data Availability Declaration}
The meshes utilized in this study will be accessible on Zenodo upon the article's acceptance. 
Given the substantial storage needs (over 6TB), 
the simulation data cannot be stored in a public repository, 
but will be provided upon reasonable request for non-commercial purposes. 
For additional inquiries, please contact the corresponding author.

\subsection*{Ethics Statement}
This study did not involve human participants or animal experiments. All investigations were conducted using a computational model of conduction system pacing capture in an interventricular septum–left bundle branch mesh.

\subsection*{Competing Interests}
None

\clearpage
\printglossaries

% ============================================== References =================================================
\clearpage

\clearpage

\bibliographystyle{unsrtnat}  % style for natbib
\bibliography{main}            % your .bib file

\begin{thebibliography}{81}
\providecommand{\natexlab}[1]{#1}
\providecommand{\url}[1]{\texttt{#1}}
\expandafter\ifx\csname urlstyle\endcsname\relax
  \providecommand{\doi}[1]{doi: #1}\else
  \providecommand{\doi}{doi: \begingroup \urlstyle{rm}\Url}\fi

\bibitem[Ghio(2004)]{ghio_interventricular_2004}
S~Ghio.
\newblock Interventricular and intraventricular dyssynchrony are common in
  heart failure patients, regardless of {QRS} duration.
\newblock \emph{European Heart Journal}, 25\penalty0 (7):\penalty0 571--578,
  April 2004.
\newblock ISSN 0195668X.
\newblock \doi{10.1016/j.ehj.2003.09.030}.
\newblock URL
  \url{https://academic.oup.com/eurheartj/article-lookup/doi/10.1016/j.ehj.2003.09.030}.

\bibitem[Auricchio et~al.(2002)Auricchio, Stellbrink, Sack, Block, Vogt,
  Bakker, Huth, Schöndube, Wolfhard, Böcker, Krahnefeld, and
  Kirkels]{auricchio_long-term_2002}
Angelo Auricchio, Christoph Stellbrink, Stefan Sack, Michael Block, J.ürgen
  Vogt, Patricia Bakker, Christof Huth, Friedrich Schöndube, Ulrich Wolfhard,
  Dirk Böcker, Olaf Krahnefeld, and Hans Kirkels.
\newblock long-term clinical effect of hemodynamically optimized cardiac
  resynchronization therapy in patients with heart failure and ventricular
  conduction delay.
\newblock \emph{Journal of the American College of Cardiology}, 39\penalty0
  (12):\penalty0 2026--2033, June 2002.
\newblock ISSN 07351097.
\newblock \doi{10.1016/S0735-1097(02)01895-8}.
\newblock URL
  \url{https://linkinghub.elsevier.com/retrieve/pii/S0735109702018958}.

\bibitem[Moss et~al.(2009)Moss, Hall, Cannom, Klein, Brown, Daubert, Estes,
  Foster, Greenberg, Higgins, Pfeffer, Solomon, Wilber, and
  Zareba]{moss_cardiac-resynchronization_2009}
Arthur~J. Moss, W.~Jackson Hall, David~S. Cannom, Helmut Klein, Mary~W. Brown,
  James~P. Daubert, N.A.~Mark Estes, Elyse Foster, Henry Greenberg, Steven~L.
  Higgins, Marc~A. Pfeffer, Scott~D. Solomon, David Wilber, and Wojciech
  Zareba.
\newblock Cardiac-{Resynchronization} {Therapy} for the {Prevention} of
  {Heart}-{Failure} {Events}.
\newblock \emph{New England Journal of Medicine}, 361\penalty0 (14):\penalty0
  1329--1338, October 2009.
\newblock ISSN 0028-4793, 1533-4406.
\newblock \doi{10.1056/NEJMoa0906431}.
\newblock URL \url{http://www.nejm.org/doi/abs/10.1056/NEJMoa0906431}.

\bibitem[Epstein et~al.(2008)Epstein, DiMarco, Ellenbogen, Estes, Freedman,
  Gettes, Gillinov, Gregoratos, Hammill, Hayes, Hlatky, Newby, Page,
  Schoenfeld, Silka, Stevenson, and Sweeney]{epstein_accahahrs_2008}
Andrew~E. Epstein, John~P. DiMarco, Kenneth~A. Ellenbogen, N.A.~Mark Estes,
  Roger~A. Freedman, Leonard~S. Gettes, A.~Marc Gillinov, Gabriel Gregoratos,
  Stephen~C. Hammill, David~L. Hayes, Mark~A. Hlatky, L.~Kristin Newby,
  Richard~L. Page, Mark~H. Schoenfeld, Michael~J. Silka, Lynne~Warner
  Stevenson, and Michael~O. Sweeney.
\newblock {ACC}/{AHA}/{HRS} 2008 {Guidelines} for {Device}-{Based} {Therapy} of
  {Cardiac} {Rhythm} {Abnormalities}.
\newblock \emph{Journal of the American College of Cardiology}, 51\penalty0
  (21):\penalty0 e1--e62, May 2008.
\newblock ISSN 07351097.
\newblock \doi{10.1016/j.jacc.2008.02.032}.
\newblock URL
  \url{https://linkinghub.elsevier.com/retrieve/pii/S0735109708007122}.

\bibitem[Vernooy et~al.(2014)Vernooy, Van~Deursen, Strik, and
  Prinzen]{vernooy_strategies_2014}
Kevin Vernooy, Caroline J.~M. Van~Deursen, Marc Strik, and Frits~W. Prinzen.
\newblock Strategies to improve cardiac resynchronization therapy.
\newblock \emph{Nature Reviews Cardiology}, 11\penalty0 (8):\penalty0 481--493,
  August 2014.
\newblock ISSN 1759-5002, 1759-5010.
\newblock \doi{10.1038/nrcardio.2014.67}.
\newblock URL \url{https://www.nature.com/articles/nrcardio.2014.67}.

\bibitem[Van Der~Wall(2014)]{van_der_wall_improvement_2014}
E.~E. Van Der~Wall.
\newblock Improvement in {CRT}: new strategies, better choices.
\newblock \emph{Netherlands Heart Journal}, 22\penalty0 (10):\penalty0
  413--414, October 2014.
\newblock ISSN 1568-5888, 1876-6250.
\newblock \doi{10.1007/s12471-014-0589-x}.
\newblock URL \url{http://link.springer.com/10.1007/s12471-014-0589-x}.

\bibitem[Ruschitzka et~al.(2013)Ruschitzka, Abraham, Singh, Bax, Borer,
  Brugada, Dickstein, Ford, Gorcsan, Gras, Krum, Sogaard, and
  Holzmeister]{ruschitzka_cardiac-resynchronization_2013}
Frank Ruschitzka, William~T. Abraham, Jagmeet~P. Singh, Jeroen~J. Bax,
  Jeffrey~S. Borer, Josep Brugada, Kenneth Dickstein, Ian Ford, John Gorcsan,
  Daniel Gras, Henry Krum, Peter Sogaard, and Johannes Holzmeister.
\newblock Cardiac-{Resynchronization} {Therapy} in {Heart} {Failure} with a
  {Narrow} {QRS} {Complex}.
\newblock \emph{New England Journal of Medicine}, 369\penalty0 (15):\penalty0
  1395--1405, October 2013.
\newblock ISSN 0028-4793, 1533-4406.
\newblock \doi{10.1056/NEJMoa1306687}.
\newblock URL \url{http://www.nejm.org/doi/10.1056/NEJMoa1306687}.

\bibitem[Parreira et~al.(2023)Parreira, Tsyganov, Artyukhina, Vernooy, Tondo,
  Adragao, Ascione, Carmo, Carvalho, Egger, Ferreira, Ghossein, Holm, Kalinin,
  Malakhova, Meine, Nunes, Podolyak, Revishvili, Shapieva, Stepanova,
  Van~Stipdonk, Taymasova, Wouters, Zubarev, Leyva, Auricchio, and
  Varma]{parreira_non-invasive_2023}
Leonor Parreira, Alexey Tsyganov, Elena Artyukhina, Kevin Vernooy, Claudio
  Tondo, Pedro Adragao, Ciro Ascione, Pedro Carmo, Salomé Carvalho, Matthias
  Egger, Antonio Ferreira, Mohammed Ghossein, Magnus Holm, Vitaly Kalinin,
  Maria Malakhova, Mathias Meine, Silvia Nunes, Dmitry Podolyak, Amiran
  Revishvili, Albina Shapieva, Vera Stepanova, Antonius Van~Stipdonk, Irina
  Taymasova, Philippe Wouters, Stepan Zubarev, Francisco Leyva, Angelo
  Auricchio, and Niraj Varma.
\newblock Non-invasive three-dimensional electrical activation mapping to
  predict cardiac resynchronization therapy response: site of latest left
  ventricular activation relative to pacing site.
\newblock \emph{EP Europace}, 25\penalty0 (4):\penalty0 1458--1466, April 2023.
\newblock ISSN 1099-5129, 1532-2092.
\newblock \doi{10.1093/europace/euad041}.
\newblock URL
  \url{https://academic.oup.com/europace/article/25/4/1458/7066898}.

\bibitem[Zweerink et~al.(2019)Zweerink, Salden, Van~Everdingen, De~Roest, Van
  De~Ven, Cramer, Doevendans, Van~Rossum, Vernooy, Prinzen, Meine, and
  Allaart]{zweerink_hemodynamic_2019}
Alwin Zweerink, Odette~A.E. Salden, Wouter~M. Van~Everdingen, Gerben~J.
  De~Roest, Peter~M. Van De~Ven, Maarten~J. Cramer, Pieter~A. Doevendans,
  Albert~C. Van~Rossum, Kevin Vernooy, Frits~W. Prinzen, Mathias Meine, and
  Cornelis~P. Allaart.
\newblock Hemodynamic {Optimization} in {Cardiac} {Resynchronization}
  {Therapy}.
\newblock \emph{JACC: Clinical Electrophysiology}, 5\penalty0 (9):\penalty0
  1013--1025, September 2019.
\newblock ISSN 2405500X.
\newblock \doi{10.1016/j.jacep.2019.05.020}.
\newblock URL
  \url{https://linkinghub.elsevier.com/retrieve/pii/S2405500X19303949}.

\bibitem[Varma et~al.(2019)Varma, Boehmer, Bhargava, Yoo, Leonelli, Costanzo,
  Saxena, Sun, Gold, Singh, Gill, and Auricchio]{varma_evaluation_2019}
Niraj Varma, John Boehmer, Kartikeya Bhargava, Dale Yoo, Fabio Leonelli,
  Mariarosa Costanzo, Anil Saxena, Lixian Sun, Michael~R. Gold, Jagmeet Singh,
  John Gill, and Angelo Auricchio.
\newblock Evaluation, {Management}, and {Outcomes} of {Patients} {Poorly}
  {Responsive} to {Cardiac} {Resynchronization} {Device} {Therapy}.
\newblock \emph{Journal of the American College of Cardiology}, 74\penalty0
  (21):\penalty0 2588--2603, November 2019.
\newblock ISSN 07351097.
\newblock \doi{10.1016/j.jacc.2019.09.043}.
\newblock URL
  \url{https://linkinghub.elsevier.com/retrieve/pii/S0735109719378039}.

\bibitem[Daubert et~al.(2016)Daubert, Behar, Martins, Mabo, and
  Leclercq]{daubert_avoiding_2016}
Claude Daubert, Nathalie Behar, Raphaël~P. Martins, Philippe Mabo, and
  Christophe Leclercq.
\newblock Avoiding non-responders to cardiac resynchronization therapy: a
  practical guide.
\newblock \emph{European Heart Journal}, page ehw270, July 2016.
\newblock ISSN 0195-668X, 1522-9645.
\newblock \doi{10.1093/eurheartj/ehw270}.
\newblock URL
  \url{http://eurheartj.oxfordjournals.org/lookup/doi/10.1093/eurheartj/ehw270}.

\bibitem[Deshmukh et~al.(2000)Deshmukh, Casavant, Romanyshyn, and
  Anderson]{deshmukh_permanent_2000}
Pramod Deshmukh, David~A. Casavant, Mary Romanyshyn, and Kathleen Anderson.
\newblock Permanent, {Direct} {His}-{Bundle} {Pacing}: {A} {Novel} {Approach}
  to {Cardiac} {Pacing} in {Patients} {With} {Normal} {His}-{Purkinje}
  {Activation}.
\newblock \emph{Circulation}, 101\penalty0 (8):\penalty0 869--877, February
  2000.
\newblock ISSN 0009-7322, 1524-4539.
\newblock \doi{10.1161/01.CIR.101.8.869}.
\newblock URL \url{https://www.ahajournals.org/doi/10.1161/01.CIR.101.8.869}.

\bibitem[Burri et~al.(2023)Burri, Jastrzebski, Cano, Čurila, de~Pooter, Huang,
  Israel, Joza, Romero, Vernooy, Vijayaraman, Whinnett, and
  Zanon]{burri_ehra_2023}
Haran Burri, Marek Jastrzebski, Óscar Cano, Karol Čurila, Jan de~Pooter,
  Weijian Huang, Carsten Israel, Jacqueline Joza, Jorge Romero, Kevin Vernooy,
  Pugazhendhi Vijayaraman, Zachary Whinnett, and Francesco Zanon.
\newblock {EHRA} clinical consensus statement on conduction system pacing
  implantation: endorsed by the {Asia} {Pacific} {Heart} {Rhythm} {Society}
  ({APHRS}), {Canadian} {Heart} {Rhythm} {Society} ({CHRS}), and {Latin}
  {American} {Heart} {Rhythm} {Society} ({LAHRS}).
\newblock \emph{Europace}, 25\penalty0 (4):\penalty0 1208--1236, April 2023.
\newblock ISSN 15322092.
\newblock \doi{10.1093/europace/euad043}.
\newblock Publisher: Oxford University Press.

\bibitem[Ali et~al.(2023{\natexlab{a}})Ali, Saqi, Arnold, Miyazawa, Keene,
  Chow, Little, Peters, Kanagaratnam, Qureshi, Ng, Linton, Lefroy, Francis,
  Boon~Lim, Tanner, Muthumala, Agarwal, Shun-Shin, Cole, and
  Whinnett]{ali_left_2023}
Nadine Ali, Khulat Saqi, Ahran~D Arnold, Alejandra~A Miyazawa, Daniel Keene,
  Ji-Jian Chow, Ian Little, Nicholas~S Peters, Prapa Kanagaratnam, Norman
  Qureshi, Fu~Siong Ng, Nick W~F Linton, David~C Lefroy, Darrel~P Francis,
  Phang Boon~Lim, Mark~A Tanner, Amal Muthumala, Girija Agarwal, Matthew~J
  Shun-Shin, Graham~D Cole, and Zachary~I Whinnett.
\newblock Left bundle branch pacing with and without anodal capture: impact on
  ventricular activation pattern and acute haemodynamics.
\newblock \emph{Europace}, 25\penalty0 (10):\penalty0 euad264, October
  2023{\natexlab{a}}.
\newblock ISSN 1099-5129, 1532-2092.
\newblock \doi{10.1093/europace/euad264}.
\newblock URL
  \url{https://academic.oup.com/europace/article/doi/10.1093/europace/euad264/7303869}.

\bibitem[Zhang et~al.(2019)Zhang, Zhou, and Gold]{zhang_left_2019}
Shu Zhang, Xiaohong Zhou, and Michael~R. Gold.
\newblock Left {Bundle} {Branch} {Pacing}.
\newblock \emph{Journal of the American College of Cardiology}, 74\penalty0
  (24):\penalty0 3039--3049, December 2019.
\newblock ISSN 07351097.
\newblock \doi{10.1016/j.jacc.2019.10.039}.
\newblock URL
  \url{https://linkinghub.elsevier.com/retrieve/pii/S0735109719382683}.

\bibitem[Vijayaraman et~al.(2023)Vijayaraman, Sharma, Cano, Ponnusamy, Herweg,
  Zanon, Jastrzebski, Zou, Chelu, Vernooy, Whinnett, Nair, Molina-Lerma,
  Curila, Zalavadia, Haseeb, Dye, Vipparthy, Brunetti, Moskal, Ross,
  Van~Stipdonk, George, Qadeer, Mumtaz, Kolominsky, Zahra, Golian, Marcantoni,
  Subzposh, and Ellenbogen]{vijayaraman_comparison_2023}
Pugazhendhi Vijayaraman, Parikshit~S. Sharma, Óscar Cano, Shunmuga~Sundaram
  Ponnusamy, Bengt Herweg, Francesco Zanon, Marek Jastrzebski, Jiangang Zou,
  Mihail~G. Chelu, Kevin Vernooy, Zachary~I. Whinnett, Girish~M. Nair, Manuel
  Molina-Lerma, Karol Curila, Dipen Zalavadia, Abdul Haseeb, Cicely Dye,
  Sharath~C. Vipparthy, Ryan Brunetti, Pawel Moskal, Alexandra Ross, Antonius
  Van~Stipdonk, Jerin George, Yusuf~K. Qadeer, Mishal Mumtaz, Jeffrey
  Kolominsky, Syeda~A. Zahra, Mehrdad Golian, Lina Marcantoni, Faiz~A.
  Subzposh, and Kenneth~A. Ellenbogen.
\newblock Comparison of {Left} {Bundle} {Branch} {Area} {Pacing} and
  {Biventricular} {Pacing} in {Candidates} for {Resynchronization} {Therapy}.
\newblock \emph{Journal of the American College of Cardiology}, 82\penalty0
  (3):\penalty0 228--241, July 2023.
\newblock ISSN 07351097.
\newblock \doi{10.1016/j.jacc.2023.05.006}.
\newblock URL
  \url{https://linkinghub.elsevier.com/retrieve/pii/S0735109723055468}.

\bibitem[Sharma et~al.(2018)Sharma, Dandamudi, Herweg, Wilson, Singh,
  Naperkowski, Koneru, Ellenbogen, and Vijayaraman]{sharma_permanent_2018}
Parikshit~S. Sharma, Gopi Dandamudi, Bengt Herweg, David Wilson, Rajeev Singh,
  Angela Naperkowski, Jayanthi~N. Koneru, Kenneth~A. Ellenbogen, and
  Pugazhendhi Vijayaraman.
\newblock Permanent {His}-bundle pacing as an alternative to biventricular
  pacing for cardiac resynchronization therapy: {A} multicenter experience.
\newblock \emph{Heart Rhythm}, 15\penalty0 (3):\penalty0 413--420, March 2018.
\newblock ISSN 15475271.
\newblock \doi{10.1016/j.hrthm.2017.10.014}.
\newblock URL
  \url{https://linkinghub.elsevier.com/retrieve/pii/S1547527117312079}.

\bibitem[Upadhyay and Tung(2017)]{upadhyay_selective_2017}
Gaurav~A. Upadhyay and Roderick Tung.
\newblock Selective versus non-selective his bundle pacing for cardiac
  resynchronization therapy.
\newblock \emph{Journal of Electrocardiology}, 50\penalty0 (2):\penalty0
  191--194, March 2017.
\newblock ISSN 00220736.
\newblock \doi{10.1016/j.jelectrocard.2016.10.003}.
\newblock URL
  \url{https://linkinghub.elsevier.com/retrieve/pii/S0022073616302734}.

\bibitem[Zhang et~al.(2018)Zhang, Guo, Hou, Wang, Qian, Li, Ge, and
  Zou]{zhang_comparison_2018}
Jinlong Zhang, Jianghong Guo, Xiaofeng Hou, Yao Wang, Zhiyong Qian, Kebei Li,
  Peibing Ge, and Jiangang Zou.
\newblock Comparison of the effects of selective and non-selective {His} bundle
  pacing on cardiac electrical and mechanical synchrony.
\newblock \emph{EP Europace}, 20\penalty0 (6):\penalty0 1010--1017, June 2018.
\newblock ISSN 1099-5129, 1532-2092.
\newblock \doi{10.1093/europace/eux120}.
\newblock URL
  \url{https://academic.oup.com/europace/article/20/6/1010/3859139}.

\bibitem[Massing and James(1976)]{massing_anatomical_1976}
G~K Massing and T~N James.
\newblock Anatomical configuration of the {His} bundle and bundle branches in
  the human heart.
\newblock \emph{Circulation}, 53\penalty0 (4):\penalty0 609--621, April 1976.
\newblock ISSN 0009-7322, 1524-4539.
\newblock \doi{10.1161/01.CIR.53.4.609}.
\newblock URL \url{https://www.ahajournals.org/doi/10.1161/01.CIR.53.4.609}.

\bibitem[Stephenson et~al.(2017)Stephenson, Atkinson, Kottas, Perde,
  Jafarzadeh, Bateman, Iaizzo, Zhao, Zhang, Anderson, Jarvis, and
  Dobrzynski]{stephenson_high_2017}
Robert~S. Stephenson, Andrew Atkinson, Petros Kottas, Filip Perde, Fatemeh
  Jafarzadeh, Mike Bateman, Paul~A. Iaizzo, Jichao Zhao, Henggui Zhang,
  Robert~H. Anderson, Jonathan~C. Jarvis, and Halina Dobrzynski.
\newblock High resolution 3-{Dimensional} imaging of the human cardiac
  conduction system from microanatomy to mathematical modeling.
\newblock \emph{Scientific Reports}, 7\penalty0 (1):\penalty0 7188, August
  2017.
\newblock ISSN 2045-2322.
\newblock \doi{10.1038/s41598-017-07694-8}.
\newblock URL \url{https://www.nature.com/articles/s41598-017-07694-8}.

\bibitem[Padala et~al.(2021)Padala, Cabrera, and
  Ellenbogen]{padala_anatomy_2021}
Santosh~K. Padala, José‐Angel Cabrera, and Kenneth~A. Ellenbogen.
\newblock Anatomy of the cardiac conduction system.
\newblock \emph{Pacing and Clinical Electrophysiology}, 44\penalty0
  (1):\penalty0 15--25, January 2021.
\newblock ISSN 0147-8389, 1540-8159.
\newblock \doi{10.1111/pace.14107}.
\newblock URL \url{https://onlinelibrary.wiley.com/doi/10.1111/pace.14107}.

\bibitem[Wu et~al.(2021)Wu, Chen, Wang, Xu, Xiao, Huang, Zheng, Jiang,
  Vijayaraman, Sharma, Su, and Huang]{wu_evaluation_2021}
Shengjie Wu, Xueying Chen, Songjie Wang, Lei Xu, Fangyi Xiao, Zhouqing Huang,
  Rujie Zheng, Limeng Jiang, Pugazhendhi Vijayaraman, Parikshit~S. Sharma, Lan
  Su, and Weijian Huang.
\newblock Evaluation of the {Criteria} to {Distinguish} {Left} {Bundle}
  {Branch} {Pacing} {From} {Left} {Ventricular} {Septal} {Pacing}.
\newblock \emph{JACC: Clinical Electrophysiology}, 7\penalty0 (9):\penalty0
  1166--1177, September 2021.
\newblock ISSN 2405500X.
\newblock \doi{10.1016/j.jacep.2021.02.018}.
\newblock URL
  \url{https://linkinghub.elsevier.com/retrieve/pii/S2405500X21002024}.

\bibitem[Sun et~al.(2022)Sun, Upadhyay, and Tung]{sun_influence_2022}
Weiping Sun, Gaurav~A. Upadhyay, and Roderick Tung.
\newblock Influence of {Capture} {Selectivity} and {Left} {Intrahisian} {Block}
  on {QRS} {Characteristics} {During} {Left} {Bundle} {Branch} {Pacing}.
\newblock \emph{JACC: Clinical Electrophysiology}, 8\penalty0 (5):\penalty0
  635--647, May 2022.
\newblock ISSN 2405500X.
\newblock \doi{10.1016/j.jacep.2022.01.012}.
\newblock URL
  \url{https://linkinghub.elsevier.com/retrieve/pii/S2405500X22000329}.

\bibitem[Bayer et~al.(2012)Bayer, Blake, Plank, and
  Trayanova]{bayer_novel_2012}
J.~D. Bayer, R.~C. Blake, G.~Plank, and N.~A. Trayanova.
\newblock A {Novel} {Rule}-{Based} {Algorithm} for {Assigning} {Myocardial}
  {Fiber} {Orientation} to {Computational} {Heart} {Models}.
\newblock \emph{Annals of Biomedical Engineering}, 40\penalty0 (10):\penalty0
  2243--2254, October 2012.
\newblock ISSN 0090-6964, 1573-9686.
\newblock \doi{10.1007/s10439-012-0593-5}.
\newblock URL \url{http://link.springer.com/10.1007/s10439-012-0593-5}.

\bibitem[Streeter et~al.(1969)Streeter, Spotnitz, Patel, Ross, and
  Sonnenblick]{streeter_fiber_1969}
Daniel~D. Streeter, Henry~M. Spotnitz, Dali~P. Patel, John Ross, and Edmund~H.
  Sonnenblick.
\newblock Fiber {Orientation} in the {Canine} {Left} {Ventricle} during
  {Diastole} and {Systole}.
\newblock \emph{Circulation Research}, 24\penalty0 (3):\penalty0 339--347,
  March 1969.
\newblock ISSN 0009-7330, 1524-4571.
\newblock \doi{10.1161/01.RES.24.3.339}.
\newblock URL \url{https://www.ahajournals.org/doi/10.1161/01.RES.24.3.339}.

\bibitem[Cabrera et~al.(2020)Cabrera, Porta-Sánchez, Tung, and
  Sánchez-Quintana]{cabrera_tracking_2020}
José-Ángel Cabrera, Andreu Porta-Sánchez, Roderick Tung, and Damián
  Sánchez-Quintana.
\newblock Tracking {Down} the {Anatomy} of the {Left} {Bundle} {Branch} to
  {Optimize} {Left} {Bundle} {Branch} {Pacing}.
\newblock \emph{JACC: Case Reports}, 2\penalty0 (5):\penalty0 750--755, May
  2020.
\newblock ISSN 26660849.
\newblock \doi{10.1016/j.jaccas.2020.04.004}.
\newblock URL
  \url{https://linkinghub.elsevier.com/retrieve/pii/S2666084920303594}.

\bibitem[Medtronic(2023)]{medtronic3830}
Medtronic.
\newblock {SelectSecure™ MRI SureScan™ 3830 Lumenless Lead – Product
  Specification Sheet}.
\newblock \url{https://europe.medtronic.com/xd-en/index.html}, 2023.
\newblock UC201702030e-selectsecure-3830-spec-sheet-en-emea-11610775.

\bibitem[Vigmond et~al.(2008)Vigmond, Weber~dos Santos, Prassl, Deo, and
  Plank]{vigmond08:_solvers}
E.~J. Vigmond, R.~Weber~dos Santos, a.~J. Prassl, M.~Deo, and G.~Plank.
\newblock Solvers for the cardiac bidomain equations.
\newblock \emph{Progress in Biophysics and Molecular Biology}, 96\penalty0
  (1-3):\penalty0 3--18, 2008.
\newblock ISSN 00796107.
\newblock \doi{10.1016/j.pbiomolbio.2007.07.012}.
\newblock ISBN: 0079-6107 (Print){\textbackslash}r0079-6107 (Linking).

\bibitem[Tusscher and Panfilov()]{tusscher_cell_2006}
K~H W J~Ten Tusscher and A~V Panfilov.
\newblock Cell model for efficient simulation of wave propagation in human
  ventricular tissue under normal and pathological conditions.
\newblock 51\penalty0 (23):\penalty0 6141--6156.
\newblock ISSN 0031-9155, 1361-6560.
\newblock \doi{10.1088/0031-9155/51/23/014}.
\newblock URL
  \url{https://iopscience.iop.org/article/10.1088/0031-9155/51/23/014}.

\bibitem[Stewart et~al.(2009)Stewart, Aslanidi, Noble, Noble, Boyett, and
  Zhang]{stewart_mathematical_2009}
Philip Stewart, Oleg~V. Aslanidi, Denis Noble, Penelope~J. Noble, Mark~R.
  Boyett, and Henggui Zhang.
\newblock Mathematical models of the electrical action potential of {Purkinje}
  fibre cells.
\newblock \emph{Philosophical Transactions of the Royal Society A:
  Mathematical, Physical and Engineering Sciences}, 367\penalty0
  (1896):\penalty0 2225--2255, June 2009.
\newblock ISSN 1364-503X, 1471-2962.
\newblock \doi{10.1098/rsta.2008.0283}.
\newblock URL
  \url{https://royalsocietypublishing.org/doi/10.1098/rsta.2008.0283}.

\bibitem[Clerc(1976)]{https://doi.org/10.1113/jphysiol.1976.sp011283}
L~Clerc.
\newblock Directional differences of impulse spread in trabecular muscle from
  mammalian heart.
\newblock \emph{The Journal of Physiology}, 255\penalty0 (2):\penalty0
  335--346, 1976.
\newblock \doi{https://doi.org/10.1113/jphysiol.1976.sp011283}.
\newblock URL
  \url{https://physoc.onlinelibrary.wiley.com/doi/abs/10.1113/jphysiol.1976.sp011283}.

\bibitem[Miklav{\v{c}}i{\v{c}} et~al.(2006)Miklav{\v{c}}i{\v{c}},
  Pav{\v{s}}elj, and Hart]{miklavvcivc2006electric}
Damijan Miklav{\v{c}}i{\v{c}}, Nata{\v{s}}a Pav{\v{s}}elj, and Francis~X Hart.
\newblock Electric properties of tissues.
\newblock \emph{Wiley encyclopedia of biomedical engineering}, 2006.

\bibitem[Flynn and O'Hagan()]{flynn_measurements_1967}
D.R. Flynn and M.E. O'Hagan.
\newblock Measurements of the thermal conductivity and electrical resistivity
  of platinum from 100 to 900 c.
\newblock 71C\penalty0 (4):\penalty0 255.
\newblock ISSN 0022-4316.
\newblock \doi{10.6028/jres.071C.021}.
\newblock URL
  \url{https://nvlpubs.nist.gov/nistpubs/jres/71C/jresv71Cn4p255_A1b.pdf}.

\bibitem[Jafarzadeh et~al.()Jafarzadeh, Farzaneh, Haddadi‐Asl, and
  Jouibari]{jafarzadeh_review_2023}
Shahabaldin Jafarzadeh, Arman Farzaneh, Vahid Haddadi‐Asl, and Iman~Sahebi
  Jouibari.
\newblock A review on electrically conductive polyurethane nanocomposites: From
  principle to application.
\newblock 44\penalty0 (12):\penalty0 8266--8302.
\newblock ISSN 0272-8397, 1548-0569.
\newblock \doi{10.1002/pc.27706}.
\newblock URL
  \url{https://4spepublications.onlinelibrary.wiley.com/doi/10.1002/pc.27706}.

\bibitem[Plank et~al.()Plank, Loewe, Neic, Augustin, Huang, Gsell, Karabelas,
  Nothstein, Prassl, Sánchez, Seemann, and Vigmond]{plank_opencarp_2021}
Gernot Plank, Axel Loewe, Aurel Neic, Christoph Augustin, Yung-Lin Huang,
  Matthias~A.F. Gsell, Elias Karabelas, Mark Nothstein, Anton~J. Prassl, Jorge
  Sánchez, Gunnar Seemann, and Edward~J. Vigmond.
\newblock The {openCARP} simulation environment for cardiac electrophysiology.
\newblock 208:\penalty0 106223.
\newblock ISSN 01692607.
\newblock \doi{10.1016/j.cmpb.2021.106223}.
\newblock URL
  \url{https://linkinghub.elsevier.com/retrieve/pii/S0169260721002972}.

\bibitem[Sobie et~al.(1997)Sobie, Susil, and Tung]{sobie_1997_generalized}
EA~Sobie, RC~Susil, and L~Tung.
\newblock A generalized activating function for predicting virtual electrodes
  in cardiac tissue.
\newblock \emph{Biophysical journal}, 73\penalty0 (September):\penalty0
  1410--1423, 1997.
\newblock URL
  \url{http://www.sciencedirect.com/science/article/pii/S0006349597781736}.
\newblock ISBN: 4109557453.

\bibitem[Rattay(1986)]{rattay_analysis_1986}
Frank Rattay.
\newblock Analysis of {Models} for {External} {Stimulation} of {Axons}.
\newblock \emph{IEEE Transactions on Biomedical Engineering}, BME-33\penalty0
  (10):\penalty0 974--977, October 1986.
\newblock ISSN 1558-2531.
\newblock \doi{10.1109/TBME.1986.325670}.
\newblock URL \url{https://ieeexplore.ieee.org/document/4122186}.
\newblock Conference Name: IEEE Transactions on Biomedical Engineering.

\bibitem[De~Pooter et~al.(2022)De~Pooter, Wauters, Van~Heuverswyn, and
  Le~Polain De~Waroux]{de_pooter_guide_2022}
Jan De~Pooter, Aurelien Wauters, Frederic Van~Heuverswyn, and Jean-Benoit
  Le~Polain De~Waroux.
\newblock A {Guide} to {Left} {Bundle} {Branch} {Area} {Pacing} {Using}
  {Stylet}-{Driven} {Pacing} {Leads}.
\newblock \emph{Frontiers in Cardiovascular Medicine}, 9:\penalty0 844152,
  February 2022.
\newblock ISSN 2297-055X.
\newblock \doi{10.3389/fcvm.2022.844152}.
\newblock URL
  \url{https://www.frontiersin.org/articles/10.3389/fcvm.2022.844152/full}.

\bibitem[Shroff et~al.(2024)Shroff, Nair, Raja, Abhilash, Fiorese, Ariyaratnam,
  Abhayaratna, Sanders, Vijayaraman, and Pathak]{shroff_comparison_2024}
Jenish~P. Shroff, Anugrah Nair, Deep~Chandh Raja, Sreevilasam~P. Abhilash,
  Simon Fiorese, Jonathan~P. Ariyaratnam, Walter~P. Abhayaratna, Prashanthan
  Sanders, Pugazhendhi Vijayaraman, and Rajeev~K. Pathak.
\newblock Comparison of {Procedural} {Outcomes} of {Lumenless} {Fixed}-{Helix}
  {Versus} {Stylet}-{Driven} {Extendable}-{Helix} {Lead} {Systems} in {Left}
  {Bundle} {Branch} {Pacing}: {COMPARE} {LBBP}.
\newblock \emph{Circulation: Arrhythmia and Electrophysiology}, 17\penalty0
  (12), December 2024.
\newblock ISSN 1941-3149, 1941-3084.
\newblock \doi{10.1161/CIRCEP.124.013385}.
\newblock URL \url{https://www.ahajournals.org/doi/10.1161/CIRCEP.124.013385}.

\bibitem[Vernooy et~al.()Vernooy, Keene, Huang, and
  Vijayaraman]{vernooy_implant_2023}
Kevin Vernooy, Daniel Keene, Weijian Huang, and Pugazhendhi Vijayaraman.
\newblock Implant, assessment, and management of conduction system pacing.
\newblock 25:\penalty0 G15--G26.
\newblock ISSN 1520-765X, 1554-2815.
\newblock \doi{10.1093/eurheartjsupp/suad115}.
\newblock URL
  \url{https://academic.oup.com/eurheartjsupp/article/25/Supplement_G/G15/7394408}.

\bibitem[James and Sherf(1971)]{james_fine_1971}
Thomas~N. James and Libi Sherf.
\newblock Fine {Structure} of the {His} {Bundle}.
\newblock \emph{Circulation}, 44\penalty0 (1):\penalty0 9--28, July 1971.
\newblock ISSN 0009-7322, 1524-4539.
\newblock \doi{10.1161/01.CIR.44.1.9}.
\newblock URL \url{https://www.ahajournals.org/doi/10.1161/01.CIR.44.1.9}.

\bibitem[Elizari(2017)]{elizari_normal_2017}
Mv~Elizari.
\newblock The normal variants in the left bundle branch system.
\newblock \emph{Journal of Electrocardiology}, 50\penalty0 (4):\penalty0
  389--399, July 2017.
\newblock ISSN 00220736.
\newblock \doi{10.1016/j.jelectrocard.2017.03.004}.
\newblock URL
  \url{https://linkinghub.elsevier.com/retrieve/pii/S0022073617300687}.

\bibitem[Titus(1973)]{titus_normal_1973}
Jack~L. Titus.
\newblock Normal {Anatomy} of the {Human} {Cardiac} {Conduction} {System}:.
\newblock \emph{Anesthesia \& Analgesia}, 52\penalty0 (4):\penalty0 508???514,
  July 1973.
\newblock ISSN 0003-2999.
\newblock \doi{10.1213/00000539-197307000-00003}.
\newblock URL \url{http://journals.lww.com/00000539-197307000-00003}.

\bibitem[Uhley and Rivkin(1959)]{uhley_visualization_1959}
Herman~N. Uhley and Laurence~M. Rivkin.
\newblock Visualization of the {Left} {Branch} of the {Human}
  {Atrioventricular} {Bundle}.
\newblock \emph{Circulation}, 20\penalty0 (3):\penalty0 419--421, September
  1959.
\newblock ISSN 0009-7322, 1524-4539.
\newblock \doi{10.1161/01.CIR.20.3.419}.
\newblock URL \url{https://www.ahajournals.org/doi/10.1161/01.CIR.20.3.419}.

\bibitem[Demoulin and Kulbertus(1972)]{demoulin_histopathological_1972}
J~C Demoulin and H~E Kulbertus.
\newblock Histopathological examination of concept of left hemiblock.
\newblock \emph{Heart}, 34\penalty0 (8):\penalty0 807--814, August 1972.
\newblock ISSN 1355-6037.
\newblock \doi{10.1136/hrt.34.8.807}.
\newblock URL \url{https://heart.bmj.com/lookup/doi/10.1136/hrt.34.8.807}.

\bibitem[Garcia‐Bustos et~al.(2017)Garcia‐Bustos, Sebastian, Izquierdo,
  Molina, Chorro, and Ruiz‐Sauri]{garciabustos_quantitative_2017}
V.~Garcia‐Bustos, R.~Sebastian, M.~Izquierdo, P.~Molina, F.~J. Chorro, and
  A.~Ruiz‐Sauri.
\newblock A quantitative structural and morphometric analysis of the {Purkinje}
  network and the {Purkinje}–myocardial junctions in pig hearts.
\newblock \emph{Journal of Anatomy}, 230\penalty0 (5):\penalty0 664--678, May
  2017.
\newblock ISSN 0021-8782, 1469-7580.
\newblock \doi{10.1111/joa.12594}.
\newblock URL \url{https://onlinelibrary.wiley.com/doi/10.1111/joa.12594}.

\bibitem[Vijayaraman et~al.(2019{\natexlab{a}})Vijayaraman, Subzposh,
  Naperkowski, Panikkath, John, Mascarenhas, Bauch, and
  Huang]{vijayaraman_prospective_2019}
Pugazhendhi Vijayaraman, Faiz~A. Subzposh, Angela Naperkowski, Ragesh
  Panikkath, Kaitlyn John, Vernon Mascarenhas, Terry~D. Bauch, and Weijian
  Huang.
\newblock Prospective evaluation of feasibility and electrophysiologic and
  echocardiographic characteristics of left bundle branch area pacing.
\newblock \emph{Heart Rhythm}, 16\penalty0 (12):\penalty0 1774--1782, December
  2019{\natexlab{a}}.
\newblock ISSN 15475271.
\newblock \doi{10.1016/j.hrthm.2019.05.011}.
\newblock URL
  \url{https://linkinghub.elsevier.com/retrieve/pii/S1547527119304424}.

\bibitem[Vijayaraman(2019)]{vijayaraman_his-bundle_2019}
Pugazhendhi Vijayaraman.
\newblock His-bundle {Pacing} to {Left} {Bundle} {Branch} {Pacing}: {Evolution}
  of {His}-{Purkinje} {Conduction} {System} {Pacing}.
\newblock \emph{Journal of Innovations in Cardiac Rhythm Management},
  0\penalty0 (5):\penalty0 3668--3673, May 2019.
\newblock ISSN 21563977, 21563993.
\newblock \doi{10.19102/icrm.2019.100504}.
\newblock URL
  \url{http://www.innovationsincrm.com/cardiac-rhythm-management/articles-2019/may/1433-his-bundle-pacing-to-left-bundle-branch-pacing}.

\bibitem[Liu et~al.(2021)Liu, Wang, Sun, Qin, and Zheng]{liu_left_2021}
Peng Liu, Qiaozhu Wang, Hongke Sun, Xinghua Qin, and Qiangsun Zheng.
\newblock Left {Bundle} {Branch} {Pacing}: {Current} {Knowledge} and {Future}
  {Prospects}.
\newblock \emph{Frontiers in Cardiovascular Medicine}, 8:\penalty0 630399,
  March 2021.
\newblock ISSN 2297-055X.
\newblock \doi{10.3389/fcvm.2021.630399}.
\newblock URL
  \url{https://www.frontiersin.org/articles/10.3389/fcvm.2021.630399/full}.

\bibitem[Ponnusamy and Vijayaraman(2021{\natexlab{a}})]{ponnusamy_how_2021}
Shunmuga~Sundaram Ponnusamy and Pugazhendhi Vijayaraman.
\newblock How to {Implant} {His} {Bundle} and {Left} {Bundle} {Pacing} {Leads}:
  {Tips} and {Pearls}.
\newblock \emph{Cardiac Failure Review}, 7:\penalty0 e13, August
  2021{\natexlab{a}}.
\newblock ISSN 20577559, 20577540.
\newblock \doi{10.15420/cfr.2021.04}.
\newblock URL \url{https://www.cfrjournal.com/articleindex/cfr.2021.04}.

\bibitem[Ponnusamy et~al.(2022)Ponnusamy, Basil, and
  Vijayaraman]{ponnusamy_electrophysiological_2022}
Shunmuga~Sundaram Ponnusamy, William Basil, and Pugazhendhi Vijayaraman.
\newblock Electrophysiological characteristics of septal perforation during
  left bundle branch pacing.
\newblock \emph{Heart Rhythm}, 19\penalty0 (5):\penalty0 728--734, May 2022.
\newblock ISSN 15475271.
\newblock \doi{10.1016/j.hrthm.2022.01.018}.
\newblock URL
  \url{https://linkinghub.elsevier.com/retrieve/pii/S1547527122000339}.

\bibitem[Ravi et~al.(2020)Ravi, Larsen, Ooms, Trohman, and
  Sharma]{ravi_late-onset_2020}
Venkatesh Ravi, Timothy Larsen, Sara Ooms, Richard Trohman, and Parikshit~S.
  Sharma.
\newblock Late-onset interventricular septal perforation from left bundle
  branch pacing.
\newblock \emph{HeartRhythm Case Reports}, 6\penalty0 (9):\penalty0 627--631,
  September 2020.
\newblock ISSN 22140271.
\newblock \doi{10.1016/j.hrcr.2020.06.008}.
\newblock URL
  \url{https://linkinghub.elsevier.com/retrieve/pii/S2214027120301214}.

\bibitem[Ghosh et~al.(2024)Ghosh, Sekar, Sriram, Sivakumar, Upadhyay, and
  Pandurangi]{ghosh_septal_2024}
Anindya Ghosh, Anbarasan Sekar, Chenni~S Sriram, Kothandam Sivakumar, Gaurav~A
  Upadhyay, and Ulhas~M Pandurangi.
\newblock Septal venous channel perforation during left bundle branch area
  pacing: a prospective study.
\newblock \emph{Europace}, 26\penalty0 (6):\penalty0 euae124, June 2024.
\newblock ISSN 1099-5129, 1532-2092.
\newblock \doi{10.1093/europace/euae124}.
\newblock URL
  \url{https://academic.oup.com/europace/article/doi/10.1093/europace/euae124/7664631}.

\bibitem[Vigmond et~al.(2021)Vigmond, Neic, Blauer, Swenson, and
  Plank]{vigmond_how_2021}
Edward~J. Vigmond, Aurel Neic, Joshua Blauer, Darrell Swenson, and Gernot
  Plank.
\newblock How {Electrode} {Position} {Affects} {Selective} {His} {Bundle}
  {Capture}: {A} {Modelling} {Study}.
\newblock \emph{IEEE Transactions on Biomedical Engineering}, 68\penalty0
  (11):\penalty0 3410--3416, November 2021.
\newblock ISSN 0018-9294, 1558-2531.
\newblock \doi{10.1109/TBME.2021.3072334}.
\newblock URL \url{https://ieeexplore.ieee.org/document/9399838/}.

\bibitem[Vijayaraman et~al.(2019{\natexlab{b}})Vijayaraman, Panikkath,
  Mascarenhas, and Bauch]{vijayaraman_left_2019}
Pugazhendhi Vijayaraman, Ragesh Panikkath, Vernon Mascarenhas, and Terry~D.
  Bauch.
\newblock Left bundle branch pacing utilizing three dimensional mapping.
\newblock \emph{Journal of Cardiovascular Electrophysiology}, 30\penalty0
  (12):\penalty0 3050--3056, December 2019{\natexlab{b}}.
\newblock ISSN 1045-3873, 1540-8167.
\newblock \doi{10.1111/jce.14242}.
\newblock URL \url{https://onlinelibrary.wiley.com/doi/10.1111/jce.14242}.

\bibitem[Orlov et~al.(2023)Orlov, Nikolaychuk, Koulouridis, Goldman, Natan,
  Armstrong, Bhattacharya, Hicks, King, and Wylie]{orlov_left_2023}
Michael~V. Orlov, Marianna Nikolaychuk, Ioannis Koulouridis, Alena Goldman,
  Shaw Natan, James Armstrong, Adhiraj Bhattacharya, Amy Hicks, Michael King,
  and John Wylie.
\newblock Left bundle area pacing: {Guiding} implant depth by ring
  measurements.
\newblock \emph{Heart Rhythm}, 20\penalty0 (1):\penalty0 55--60, January 2023.
\newblock ISSN 15475271.
\newblock \doi{10.1016/j.hrthm.2022.09.013}.
\newblock URL
  \url{https://linkinghub.elsevier.com/retrieve/pii/S1547527122024195}.

\bibitem[Ponnusamy and Vijayaraman(2021{\natexlab{b}})]{ponnusamy_late_2021}
Shunmuga~Sundaram Ponnusamy and Pugazhendhi Vijayaraman.
\newblock Late dislodgement of left bundle branch pacing lead and successful
  extraction.
\newblock \emph{Journal of Cardiovascular Electrophysiology}, 32\penalty0
  (8):\penalty0 2346--2349, August 2021{\natexlab{b}}.
\newblock ISSN 1045-3873, 1540-8167.
\newblock \doi{10.1111/jce.15155}.
\newblock URL \url{https://onlinelibrary.wiley.com/doi/10.1111/jce.15155}.

\bibitem[Moustafa et~al.(2023)Moustafa, Tang, and
  Khan]{moustafa_conduction_2023}
Ahmed~T. Moustafa, Anthony~Sl. Tang, and Habib~Rehman Khan.
\newblock Conduction system pacing on track to replace {CRT}? {Review} of
  current evidence and prospects of conduction system pacing.
\newblock \emph{Frontiers in Cardiovascular Medicine}, 10:\penalty0 1220709,
  August 2023.
\newblock ISSN 2297-055X.
\newblock \doi{10.3389/fcvm.2023.1220709}.
\newblock URL
  \url{https://www.frontiersin.org/articles/10.3389/fcvm.2023.1220709/full}.

\bibitem[Glikson et~al.(2025)Glikson, Burri, Abdin, Cano, Curila, De~Pooter,
  Diaz, Drossart, Huang, Israel, Jastrzębski, Joza, Karvonen, Keene, Leclercq,
  Mullens, Pujol-Lopez, Rao, Vernooy, Vijayaraman, Zanon, Michowitz, Nielsen,
  Boersma, Blomström-Lundqvist, Kronborg, Chung, Tse, Khan, Leyva,
  Rojel-Martinez, Ruciński, and Varma]{glikson_european_2025}
Michael Glikson, Haran Burri, Amr Abdin, Oscar Cano, Karol Curila, Jan
  De~Pooter, Juan~C Diaz, Inga Drossart, Weijian Huang, Carsten~W Israel, Marek
  Jastrzębski, Jacqueline Joza, Jarkko Karvonen, Daniel Keene, Christophe
  Leclercq, Wilfried Mullens, Margarida Pujol-Lopez, Archana Rao, Kevin
  Vernooy, Pugazhendhi Vijayaraman, Francesco Zanon, Yoav Michowitz,
  Jens~Cosedis Nielsen, Lucas Boersma, Carina Blomström-Lundqvist, Mads~Brix
  Kronborg, Mina~K Chung, Hung~Fat Tse, Habib~Rehman Khan, Francisco Leyva,
  Ulises Rojel-Martinez, Marcin Ruciński, and Niraj Varma.
\newblock European {Society} of {Cardiology} ({ESC}) clinical consensus
  statement on indications for conduction system pacing, with special
  contribution of the {European} {Heart} {Rhythm} {Association} of the {ESC}
  and endorsed by the {Asia} {Pacific} {Heart} {Rhythm} {Society}, the
  {Canadian} {Heart} {Rhythm} {Society}, the {Heart} {Rhythm} {Society}, and
  the {Latin} {American} {Heart} {Rhythm} {Society}.
\newblock \emph{Europace}, 27\penalty0 (4):\penalty0 euaf050, March 2025.
\newblock ISSN 1099-5129, 1532-2092.
\newblock \doi{10.1093/europace/euaf050}.
\newblock URL
  \url{https://academic.oup.com/europace/article/doi/10.1093/europace/euaf050/8100402}.

\bibitem[Heckman et~al.()Heckman, Luermans, Salden, van Stipdonk, Mafi-Rad,
  Prinzen, and Vernooy]{heckman_physiology_2021}
Luuk Heckman, Justin Luermans, Floor Salden, Antonius Martinus~Wilhelmus van
  Stipdonk, Masih Mafi-Rad, Frits Prinzen, and Kevin Vernooy.
\newblock Physiology and practicality of left ventricular septal pacing.
\newblock 10\penalty0 (3):\penalty0 165--171.
\newblock ISSN 2050-3369.
\newblock \doi{10.15420/aer.2021.21}.

\bibitem[Curila et~al.()Curila, Poviser, Stros, Jurak, Whinnett, Jastrzebski,
  Waldauf, Smisek, Viscor, Hozman, Osmancik, Kryze, and
  Kautzner]{curila_lvsp_2024}
Karol Curila, Lukas Poviser, Petr Stros, Pavel Jurak, Zachary Whinnett, Marek
  Jastrzebski, Petr Waldauf, Radovan Smisek, Ivo Viscor, Marek Hozman, Pavel
  Osmancik, Lukas Kryze, and Josef Kautzner.
\newblock {LVSP} and {LBBP} result in similar or improved {LV} synchrony and
  hemodynamics compared to {BVP}.
\newblock 10\penalty0 (7):\penalty0 1722--1732.
\newblock ISSN 2405500X.
\newblock \doi{10.1016/j.jacep.2024.04.022}.
\newblock URL
  \url{https://linkinghub.elsevier.com/retrieve/pii/S2405500X24003487}.

\bibitem[Wang et~al.(2020)Wang, Gu, Qian, Hou, Chen, Qiu, Jiang, Zhang, Wu,
  Chen, and Zou]{wang_efficacy_2020}
Yao Wang, Kai Gu, Zhiyong Qian, Xiaofeng Hou, Xing Chen, Yuanhao Qiu, Zeyu
  Jiang, Xinwei Zhang, Hongping Wu, Minglong Chen, and Jiangang Zou.
\newblock The efficacy of left bundle branch area pacing compared with
  biventricular pacing in patients with heart failure: {A} matched
  case–control study.
\newblock \emph{Journal of Cardiovascular Electrophysiology}, 31\penalty0
  (8):\penalty0 2068--2077, August 2020.
\newblock ISSN 1045-3873, 1540-8167.
\newblock \doi{10.1111/jce.14628}.
\newblock URL \url{https://onlinelibrary.wiley.com/doi/10.1111/jce.14628}.

\bibitem[Cano and Vijayaraman(2021)]{cano_left_2021}
Óscar Cano and Pugazhendhi Vijayaraman.
\newblock Left {Bundle} {Branch} {Area} {Pacing}: {Implant} {Technique},
  {Definitions}, {Outcomes}, and {Complications}.
\newblock \emph{Current Cardiology Reports}, 23\penalty0 (11):\penalty0 155,
  November 2021.
\newblock ISSN 1523-3782, 1534-3170.
\newblock \doi{10.1007/s11886-021-01585-1}.
\newblock URL \url{https://link.springer.com/10.1007/s11886-021-01585-1}.

\bibitem[Hopman et~al.(2023)Hopman, Beunder, Borodzicz-Jazdzyk, Götte, and
  Van~Halm]{hopman_loss_2023}
Luuk~H.G.A. Hopman, Kyle~P. Beunder, Sonia Borodzicz-Jazdzyk, Marco~J.W.
  Götte, and Vokko~P. Van~Halm.
\newblock Loss of capture of conduction system pacemaker caused by fibrosis
  surrounding the lead: a case report.
\newblock \emph{BMC Cardiovascular Disorders}, 23\penalty0 (1):\penalty0 621,
  December 2023.
\newblock ISSN 1471-2261.
\newblock \doi{10.1186/s12872-023-03656-3}.
\newblock URL
  \url{https://bmccardiovascdisord.biomedcentral.com/articles/10.1186/s12872-023-03656-3}.

\bibitem[Arnold et~al.(2020)Arnold, Whinnett, and
  Vijayaraman]{arnold_his-purkinje_2020}
Ahran~D. Arnold, Zachary~I. Whinnett, and Pugazhendhi Vijayaraman.
\newblock His-{Purkinje} conduction system pacing: {State} of the art in 2020.
\newblock \emph{Arrhythmia and Electrophysiology Review}, 9\penalty0 (3), 2020.
\newblock ISSN 20503377.
\newblock \doi{10.15420/AER.2020.14}.
\newblock Publisher: Radcliffe Cardiology.

\bibitem[Tamborero et~al.(2006)Tamborero, Mont, Alanis, Berruezo, Tolosana,
  Sitges, Vidal, and Brugada]{tamborero_anodal_2006}
David Tamborero, Lluis Mont, Roberto Alanis, Antonio Berruezo, Jose~Maria
  Tolosana, Marta Sitges, Barbara Vidal, and Josep Brugada.
\newblock Anodal {Capture} in {Cardiac} {Resynchronization} {Therapy}
  {Implications} for {Device} {Programming}.
\newblock \emph{Pacing and Clinical Electrophysiology}, 29\penalty0
  (9):\penalty0 940--945, September 2006.
\newblock ISSN 0147-8389, 1540-8159.
\newblock \doi{10.1111/j.1540-8159.2006.00466.x}.
\newblock URL
  \url{https://onlinelibrary.wiley.com/doi/10.1111/j.1540-8159.2006.00466.x}.

\bibitem[Wikswo et~al.(1995)Wikswo, Lin, and Abbas]{wikswo_virtual_1995}
J.P. Wikswo, S.F. Lin, and R.A. Abbas.
\newblock Virtual electrodes in cardiac tissue: a common mechanism for anodal
  and cathodal stimulation.
\newblock \emph{Biophysical Journal}, 69\penalty0 (6):\penalty0 2195--2210,
  December 1995.
\newblock ISSN 00063495.
\newblock \doi{10.1016/S0006-3495(95)80115-3}.
\newblock URL
  \url{https://linkinghub.elsevier.com/retrieve/pii/S0006349595801153}.

\bibitem[Sambelashvili et~al.()Sambelashvili, Nikolski, and
  Efimov]{sambelashvili_virtual_2004}
Aleksandre~T. Sambelashvili, Vladimir~P. Nikolski, and Igor~R. Efimov.
\newblock Virtual electrode theory explains pacing threshold increase caused by
  cardiac tissue damage.
\newblock 286\penalty0 (6):\penalty0 H2183--H2194.
\newblock ISSN 0363-6135, 1522-1539.
\newblock \doi{10.1152/ajpheart.00637.2003}.
\newblock URL \url{https://www.physiology.org/doi/10.1152/ajpheart.00637.2003}.

\bibitem[Vijayaraman et~al.(2022)Vijayaraman, Cano, Ponnusamy, Molina-Lerma,
  Chan, Padala, Sharma, Whinnett, Herweg, Upadhyay, Subzposh, Patel, Beer,
  Bednarek, Kielbasa, Tung, Ellenbogen, and Jastrzebski]{vijayaraman_left_2022}
Pugazhendhi Vijayaraman, Oscar Cano, Shunmuga~Sundaram Ponnusamy, Manuel
  Molina-Lerma, Joseph~Y.S. Chan, Santosh~K. Padala, Parikshit~S. Sharma,
  Zachary~I. Whinnett, Bengt Herweg, Gaurav~A. Upadhyay, Faiz~A. Subzposh,
  Neil~R. Patel, Dominik~A. Beer, Agnieszka Bednarek, Grzegorz Kielbasa,
  Roderick Tung, Kenneth~A. Ellenbogen, and Marek Jastrzebski.
\newblock Left bundle branch area pacing in patients with heart failure and
  right bundle branch block: {Results} from {International} {LBBAP}
  {Collaborative}-{Study} {Group}.
\newblock \emph{Heart Rhythm O2}, 3\penalty0 (4):\penalty0 358--367, August
  2022.
\newblock ISSN 26665018.
\newblock \doi{10.1016/j.hroo.2022.05.004}.
\newblock URL
  \url{https://linkinghub.elsevier.com/retrieve/pii/S266650182200109X}.

\bibitem[Huang et~al.(2017)Huang, Su, Wu, Xu, Xiao, Zhou, and
  Ellenbogen]{huang_novel_2017}
Weijian Huang, Lan Su, Shengjie Wu, Lei Xu, Fangyi Xiao, Xiaohong Zhou, and
  Kenneth~A. Ellenbogen.
\newblock A {Novel} {Pacing} {Strategy} {With} {Low} and {Stable} {Output}:
  {Pacing} the {Left} {Bundle} {Branch} {Immediately} {Beyond} the {Conduction}
  {Block}.
\newblock \emph{Canadian Journal of Cardiology}, 33\penalty0 (12):\penalty0
  1736.e1--1736.e3, December 2017.
\newblock ISSN 0828282X.
\newblock \doi{10.1016/j.cjca.2017.09.013}.
\newblock URL
  \url{https://linkinghub.elsevier.com/retrieve/pii/S0828282X17310322}.

\bibitem[Vinther et~al.(2024)Vinther, Sandgaard, Risum, and
  Philbert]{vinther_late_2024}
Michael Vinther, Niels~C.F. Sandgaard, Niels Risum, and Berit~Th. Philbert.
\newblock Late perforation of a left bundle branch area pacing lead causing
  ventricular fibrillation: {A} case report.
\newblock \emph{HeartRhythm Case Reports}, 10\penalty0 (7):\penalty0 509--513,
  July 2024.
\newblock ISSN 22140271.
\newblock \doi{10.1016/j.hrcr.2024.05.002}.
\newblock URL
  \url{https://linkinghub.elsevier.com/retrieve/pii/S221402712400099X}.

\bibitem[Wang et~al.(2021)Wang, Lan, Zhang, Zheng, Gao, Bai, Wu, Xu, Wang, and
  Xu]{wang_lbbap_2021}
Shaoxian Wang, Rongfang Lan, Ning Zhang, Jia Zheng, Yuan Gao, Jian Bai, Xiang
  Wu, Xinyue Xu, Tianqi Wang, and Wei Xu.
\newblock {LBBAP} in patients with normal intrinsic {QRS} duration:
  {Electrical} and mechanical characteristics.
\newblock \emph{Pacing and Clinical Electrophysiology}, 44\penalty0
  (1):\penalty0 82--92, January 2021.
\newblock ISSN 0147-8389, 1540-8159.
\newblock \doi{10.1111/pace.14114}.
\newblock URL \url{https://onlinelibrary.wiley.com/doi/10.1111/pace.14114}.

\bibitem[Huang et~al.(2019)Huang, Chen, Su, Wu, Xia, and
  Vijayaraman]{huang_beginners_2019}
Weijian Huang, Xueying Chen, Lan Su, Shengjie Wu, Xue Xia, and Pugazhendhi
  Vijayaraman.
\newblock A beginner's guide to permanent left bundle branch pacing.
\newblock \emph{Heart Rhythm}, 16\penalty0 (12):\penalty0 1791--1796, December
  2019.
\newblock ISSN 15475271.
\newblock \doi{10.1016/j.hrthm.2019.06.016}.
\newblock URL
  \url{https://linkinghub.elsevier.com/retrieve/pii/S1547527119305739}.

\bibitem[Vijayaraman(2020)]{vijayaraman_deep_2020}
Pugazhendhi Vijayaraman.
\newblock Deep septal, distal {His} bundle pacing for cardiac resynchronization
  therapy.
\newblock \emph{HeartRhythm Case Reports}, 6\penalty0 (10):\penalty0 791--793,
  October 2020.
\newblock ISSN 22140271.
\newblock \doi{10.1016/j.hrcr.2020.07.024}.
\newblock URL
  \url{https://linkinghub.elsevier.com/retrieve/pii/S2214027120301743}.

\bibitem[Diaz et~al.(2023)Diaz, Duque, Aristizabal, Marin, Niño, Bastidas,
  Ruiz, Matos, Hoyos, Hincapie, Velasco, and Romero]{diaz_emerging_2023}
Juan~Carlos Diaz, Mauricio Duque, Julian Aristizabal, Jorge Marin, Cesar Niño,
  Oriana Bastidas, Luis~Miguel Ruiz, Carlos~D Matos, Carolina Hoyos, Daniela
  Hincapie, Alejandro Velasco, and Jorge~E Romero.
\newblock The {Emerging} {Role} of {Left} {Bundle} {Branch} {Area} {Pacing} for
  {Cardiac} {Resynchronisation} {Therapy}.
\newblock \emph{Arrhythmia \& Electrophysiology Review}, 12:\penalty0 e29,
  December 2023.
\newblock ISSN 20503377, 20503369.
\newblock \doi{10.15420/aer.2023.15}.
\newblock URL \url{https://www.aerjournal.com/articleindex/aer.2023.15}.

\bibitem[Saxonhouse et~al.(2005)Saxonhouse, Conti, and
  Curtis]{saxonhouse_current_2005}
Sherry~J. Saxonhouse, Jamie~B. Conti, and Anne~B. Curtis.
\newblock Current of injury predicts adequate active lead fixation in permanent
  pacemaker/defibrillation leads.
\newblock \emph{Journal of the American College of Cardiology}, 45\penalty0
  (3):\penalty0 412--417, February 2005.
\newblock ISSN 07351097.
\newblock \doi{10.1016/j.jacc.2004.10.045}.
\newblock URL
  \url{https://linkinghub.elsevier.com/retrieve/pii/S0735109704021588}.

\bibitem[Shali et~al.(2022)Shali, Wu, Bai, Wang, Qin, Wang, Liang, Chen, Su,
  Chen, and Ge]{shali_current_2022}
Shalaimaiti Shali, Weiyun Wu, Jin Bai, Wei Wang, Shengmei Qin, Jingfeng Wang,
  Yixiu Liang, Haiyan Chen, Yangang Su, Xueying Chen, and Junbo Ge.
\newblock Current of injury is an indicator of lead depth and performance
  during left bundle branch pacing lead implantation.
\newblock \emph{Heart Rhythm}, 19\penalty0 (8):\penalty0 1281--1288, August
  2022.
\newblock ISSN 15475271.
\newblock \doi{10.1016/j.hrthm.2022.04.027}.
\newblock URL
  \url{https://linkinghub.elsevier.com/retrieve/pii/S1547527122019518}.

\bibitem[Su et~al.(2020)Su, Xu, Cai, Xu, Vijayaraman, Sharma, Chen, Zheng, Wu,
  and Huang]{su_electrophysiological_2020}
Lan Su, Tiancheng Xu, Mengxing Cai, Lei Xu, Pugazhendhi Vijayaraman,
  Parikshit~S. Sharma, Xiao Chen, Rujie Zheng, Shengjie Wu, and Weijian Huang.
\newblock Electrophysiological characteristics and clinical values of left
  bundle branch current of injury in left bundle branch pacing.
\newblock \emph{Journal of Cardiovascular Electrophysiology}, 31\penalty0
  (4):\penalty0 834--842, April 2020.
\newblock ISSN 1045-3873, 1540-8167.
\newblock \doi{10.1111/jce.14377}.
\newblock URL \url{https://onlinelibrary.wiley.com/doi/10.1111/jce.14377}.

\bibitem[Ali et~al.(2023{\natexlab{b}})Ali, Arnold, Miyazawa, Keene, Chow,
  Little, Peters, Kanagaratnam, Qureshi, Ng, Linton, Lefroy, Francis,
  Phang~Boon, Tanner, Muthumala, Shun-Shin, Cole, and
  Whinnett]{ali_comparison_2023}
Nadine Ali, Ahran~D Arnold, Alejandra~A Miyazawa, Daniel Keene, Ji-Jian Chow,
  Ian Little, Nicholas~S Peters, Prapa Kanagaratnam, Norman Qureshi, Fu~Siong
  Ng, Nick W~F Linton, David~C Lefroy, Darrel~P Francis, Lim Phang~Boon, Mark~A
  Tanner, Amal Muthumala, Matthew~J Shun-Shin, Graham~D Cole, and Zachary~I
  Whinnett.
\newblock Comparison of methods for delivering cardiac resynchronization
  therapy: an acute electrical and haemodynamic within-patient comparison of
  left bundle branch area, {His} bundle, and biventricular pacing.
\newblock \emph{EP Europace}, 25\penalty0 (3):\penalty0 1060--1067, March
  2023{\natexlab{b}}.
\newblock ISSN 1099-5129, 1532-2092.
\newblock \doi{10.1093/europace/euac245}.
\newblock URL
  \url{https://academic.oup.com/europace/article/25/3/1060/7025387}.

\bibitem[Roth(2002)]{roth2002electrical}
Bradley~J Roth.
\newblock Electrical conductivity values used with the bidomain model of
  cardiac tissue.
\newblock \emph{IEEE Transactions on Biomedical Engineering}, 44\penalty0
  (4):\penalty0 326--328, 2002.

\end{thebibliography}

\end{document}